\documentclass[english,a4paper]{article}
\usepackage{dcolumn}
\usepackage{bm}
\usepackage{graphicx}
\usepackage{amsmath}
\usepackage{amsfonts}
\usepackage{amssymb}
\usepackage{amsthm}
\usepackage{amstext}
\usepackage{amsbsy}
\usepackage{amsopn}
\usepackage{amscd}
\usepackage{amsxtra}
\usepackage{mathtools}
\usepackage{commath}

\def\openone{\leavevmode\hbox{\small1\kern-3.3pt\normalsize1}}

\usepackage[usenames,dvipsnames]{color}

\def\phi{\varphi}
\def\epsilon{\varepsilon}

\DeclareMathOperator{\e}{e}

\def\P{\mathtt P}

\def\a{\mathfrak a}

\begin{document}

\title{A comparison between optimal control and shortcut to adiabaticity protocols in a linear control system}

\author{V. Martikyan, D. Gu\'ery-Odelin\footnote{Laboratoire de Collisions Agr\'egats R\'eactivit\'e, Universit\'e Paul Sabatier, 118 Route de Narbonne, 31062 Toulouse Cedex 4, France}, D. Sugny\footnote{Laboratoire Interdisciplinaire Carnot de
Bourgogne (ICB), UMR 6303 CNRS-Universit\'e Bourgogne-Franche Comt\'e, 9 Av. A.
Savary, BP 47 870, F-21078 Dijon Cedex, France, dominique.sugny@u-bourgogne.fr}}

\maketitle

\begin{abstract}
We analyze the control of the motion of a charged particle by means of an external electric field. The system is constrained to move along a given direction. The goal of the control is to change the speed of the particle in a fixed time with zero initial and final accelerations, while minimizing a cost functional in order to achieve a smooth transport of the system. We solve this linear control problem by using shortcut to adiabaticity methods and optimal control techniques. STA protocols are built upon local constraints. By extending the number of space variables, we explain how optimal control protocols can accommodate for such constraints and provide robust solutions against slight initial and final time uncertainties. Conversely, optimal control can guide the choice of the class of functions involved in STA processes.  Such mutual benefices of the two approaches remain valid beyond the specific example studied here.
\end{abstract}

\section{Introduction}\label{sec1}
In physics and engineering, control theory provides a systematic way for driving a dynamical system from a given initial state into a desired target state~\cite{glaserreview,brifreview,altafinireview,dongreview,alessandrobook,RMP19}. In this framework, optimal control theory (OCT) is a general mathematical procedure which allows to design external fields or a sequence of pulses, while minimizing or maximizing specific functionals such as the control time or the used energy~\cite{glaserreview}. A rigorous setting for Optimal control theory was given with the Pontryagin Maximum Principle (PMP) in the late 1950’s~\cite{pont}. Its development was originally inspired by problems of space dynamics, but it is now a key tool to study a large spectrum of applications ranging from robotics, economics, and mechanics~\cite{bonnardbook}. More recently, different applications at the microscopic and quantum scales have been also developed~\cite{glaserreview,brifreview,altafinireview,dongreview,alessandrobook}. Two different types of approaches can be used to solve optimal control problems, namely geometric~\cite{bonnardbook,jurdjevicbook} and numerical methods for dynamical systems~\cite{grape,reichkrotov,gross,stefanatos:2009} of low and high dimensions, respectively. Geometric techniques have different interests which are complementary to numerical optimization procedures. In particular, rigorous results of optimality can be established. The global optimal control field can be designed analytically or, at least, with a very high numerical precision. The geometric methods are also interesting to unreveal the physical mechanisms used by the control process~\cite{alessandro,boscain,garon,khanejaspin,bonnard:2012,lapert:2012}. In contrast, numerical algorithms do not generally provide physical insight about the dynamics of the controlled system or the structure of the control field~\cite{glaserreview}.

Other approaches have been proposed to date to solve control problems with relatively simple fields, which can be derived analytically. They extend from intuitive schemes such as sudden or adiabatic methods~\cite{adiabaticreview,stirapRMP}, which are valid in some specific limits, to more elaborate pulse sequences based on Shortcut To Adiabaticity (STA) techniques~\cite{reviewSTA1,reviewSTA2,reviewSTA3,reviewSTA4,STA,STAnjp,daemsprl,vandamme}. The original motivation of STA protocols is to speed up adiabatic control of the dynamical system, while preserving as much as possible its efficiency and robustness. Control techniques have been largely explored in the past few years in physics. STA and OCT have been, for instance, recently applied in the context of quantum thermodynamics~(see \cite{deffnerbook,kosloff:2013,sivak:2012,acconcia:2015,cavina:2018,abah:2018} to mention a few  and references therein). However, most of the studies focus on the control of non-linear systems, and very little in the linear case~\cite{bryson,bressan,liberzon,brockett,Lithesis,li:2011}. In particular, different works have compared STA and OCT and shown how to combine them in this non-linear setting~\cite{stefanatos:2010,xen:2011,guery:2014}. We propose in this paper to revisit this comparison for a linear control system. From a mathematical point of view, such systems are particularly appealing because analytical solutions can be easily found and optimal control problems are easier to solve~\cite{bryson,bressan,liberzon,Lithesis,li:2011}. Physically, this analysis is useful because non-linear dynamics can be approximated to some extent by linear ones. In addition, it has been shown in the case of spin systems that a non-linear mapping can be established between non-linear and linear control dynamics~\cite{li:2017}. This mapping, which allows to design analytical and efficient broadband pulses for spin dynamics from the control of springs, reinforces the interest of studying linear systems. Such studies are also interesting in specific domains where model systems are linear. An example is given by two-dimensional Fourier transform ion cyclotron resonance mass spectroscopy for which control techniques can significantly advance the sensitivity and the efficiency of current excitation and detection processes~\cite{bodenhausen:2016,delsuc:2013,delsuc:2016,delsuc:2017}. In this control process, ions are subjected to a static magnetic field and to a perpendicular time-dependent electric field which is shaped to manipulate ion trajectories and to steer the system to specific final position and speed. To the best of our knowledge, OCT and STA have not been applied in this domain. As explained below, the control protocol we study in this paper can be viewed as a first step towards the control of ion trajectories.

In order to explain how both techniques can benefit mutually from each other, we consider a minimal linear model involving the control of a classical one-dimensional charged particle by means of an electric field. The goal of the control is to change the speed of the system in a given time with the constraint that the initial and final accelerations are zero. The optimal control problem is defined through a cost functional, which allows a smooth transport of the particle. The decisive advantage of this simple control scenario is that a complete geometric and analytical description can be carried out for the two approaches. In particular, different bases of functions to expand the solutions can be derived in the STA case, singular and regular optimal fields (see below for the definition) can be also determined analytically. The relative efficiency of the different control protocols is measured with respect to the global optimal solution. We show that OCT can guide the choice of the basis of functions in which the STA solution is expanded. This approach can be generalized and the same analysis can be conducted in the case of more constrained boundary conditions where the time derivatives of the speed at any order have to be zero.

The paper is organized as follows. We present the model system and the control problem in Sec.~\ref{sec2}. Different STA protocols using polynomial, trigonometric and real exponential function bases are proposed in Sec.~\ref{sec3}. Section~\ref{sec4} is dedicated to the derivation of the singular optimal solution. A regularization is proposed in Sec.~\ref{sec5} to achieve finite-amplitude control fields. A thorough comparison is made between STA and OCT pulses. A generalization in which the time derivatives of the speed at any order have to be zero at the initial and final times of the process is considered in Sec.~\ref{sec6}. Conclusion and prospective views are given in Sec.~\ref{sec7}. Technical details are reported in the Supplementary Material.
\section{The model system}\label{sec2}
We consider the control of the dynamics of a charged particle by means of an external electric field. The control process is aimed at accelerating the particle to a given speed and we assume that the system can only move along one direction. The position of the system is given by the real coordinate $X\in\mathbb{R}$. The system is also subjected to friction modeled by a force proportional to the speed of the system. Using Newton's law along the $X$- axis, it is straightforward to show that:
\begin{equation}
m\ddot{X}+\eta \dot{X}=qE(t),
\end{equation}
where $m$ and $q$ are respectively the mass and the charge of the particle, and $\eta$ the friction coefficient. The initial and target states are $\dot{X}(0)=0$ and $\dot{X}(T)=V$ where $T$ is the fixed control time and $V$ the final speed. The goal of the control is to reach exactly the target state at time $T$, while minimizing a cost functional $\mathcal{C}$ which can be expressed as the sum of the time average of the speed and acceleration along the trajectory:
$$
\mathcal{C}=\frac{1}{T}\int_0^T (\ddot{X}^2+\alpha\dot{X}^2)dt,
$$
where $\alpha$ is a positive parameter to express the relative weight between the two terms of $\mathcal{C}$. This choice ensures a smooth transport of the particle, avoiding, e.g., large speed variations. An additional constraint on the initial and final accelerations, $\ddot{X}(0)=\ddot{X}(T)=0$ is also accounted for. Using normalized coordinates, the dynamics of the system can be rewritten as:
\begin{equation}\label{eqstart}
\dot{x}+x=u(t)
\end{equation}
where $x$ is the normalized speed and $u(t)$ the normalized control field used to drive the particle from $x(0)=0$ to $x(T)=1$ at a normalized time $T$. The control has also to satisfy the constraints $\dot{x}(0)=0$ and $\dot{x}(T)=0$, while minimizing the running cost $\mathcal{C}=\int_0^T[x(t)^2+\dot{x}(t)^2]dt$. Without loss of generality, we will assume that $T=1$ in the rest of the paper.

\section{Shortcut to Adiabaticity protocols}\label{sec3}

STA protocols have originally been developed in the context of quantum mechanics~\cite{reviewSTA1,reviewSTA2,reviewSTA3,reviewSTA4}. They correspond to fast routes between initial and final states that are connected through a slow (adiabatic) time evolution when a controlling parameter is changed in time. In the example under study, an adiabatic process corresponds to a very slow time-dependent variation of the control field $u(t)$ from 0 to 1. Neglecting the time derivative $\dot{x}(t)$ in Eq.~\eqref{eqstart}, we obtain that the state of the system $x(t)$ can be identified to $u(t)$. The different STA methods exploit the algebraic structure of quantum mechanics. STA has been extended to statistical physics and classical mechanics. The simplest method used in this latter context is reminiscent whilst slightly different of motion planning method in control theory~\cite{motionplanning}, and is commonly referred to as inverse engineering~\cite{reviewSTA4}. It consists in extrapolating the trajectory from the required boundary conditions and inferring from this interpolation the shape of the control field that should be applied to the system under study.

We exemplify this approach in this paragraph by deriving a few simple analytical solutions of Eq.~\eqref{eqstart}. We consider three different ansatz for which the trajectory of the system $x(t)$ belongs to a family of functions, namely polynomial, trigonometric or exponential functions. This choice is arbitrary and the same study could be made for other families.  The trajectory depends on a few number of parameters which are determined from the boundary conditions and correspond to the minimum of the cost functional $\mathcal{C}$.

For the polynomial function family, we assume that the trajectory can be expressed as a polynomial function in the time $t$ up to a given order $N$. The constraints $x(0)=\dot{x}(0)=0$ impose that $a_0=a_1=0$, and therefore:
\begin{equation}
\label{xp}
x(t)=\sum\limits_{k=2}^N a_k t^k.
\end{equation}
From $x(T=1)=1$ and $\dot{x}(T=1)=0$, we deduce that:
\begin{equation}
\label{xp1}
\begin{cases}
\sum\limits_{k=2}^N a_k=1
\\
\sum\limits_{k=2}^N k a_k=0
\end{cases}
\end{equation}
There are two equations and $N-1$ unknowns, so $N-3$ free parameters. For instance, at the order $N=3$, all the coefficients can be determined from the boundary conditions and the unique solution $x_3(t)$ is:
$$
x_3(t)=3t^2-2t^3.
$$
At the order 4, it can be shown that a one-parameter family of polynomials $x_4(t)|_a$ is solution of the control problem:
$$
x_4(t)=(3+a)t^2-(2+2a)t^3+a t^4,
$$
where $a$ is a free real parameter, which is chosen to minimize $\mathcal{C}$. Numerical computation leads to $a=-0.8076923$. The same study can be done for polynomials of higher orders. The polynomials depend respectively on two and three free parameters at the orders 5 and 6. Technical details can be found in Supplementary Sec. I. The parameters for $N$ going from 4 to 6 are given in Tab.~\ref{tab1}.

A trigonometric expansion can also be used to express the trajectory $x(t)$:
\begin{equation}
\label{xtg}
x(t)=\sum\limits_{k=1}^N a_k \sin\left(\frac{k\pi}{2} t\right).
\end{equation}
Note that $x(0)=0$ is automatically satisfied, and the condition $\dot{x}(0)=0$ leads to:
\begin{eqnarray}
\sum\limits_{k=1}^N k a_k=0.
\end{eqnarray}
From $x(1)=1$ and $\dot{x}(1)=0$, we get:
\begin{equation}
\label{xtg1}
\begin{cases}
\sum\limits_{k=1}^N a_k \sin\left(\frac{k\pi}{2}\right)=1,
\\
\sum\limits_{k=1}^N k a_k \cos\left(\frac{k\pi}{2}\right)=0.
\end{cases}
\end{equation}
We have three equations and $N$ unknowns, hence there are again $N-3$ free parameters to find (see Tab.~\ref{tab1} for the solutions in the cases $N=4$, 5 and 6).
\begin{table}
\centering
\begin{tabular}{|c|c|c|c|c|}
                 \hline
                 & $N$ 	 & $a$ & $b$ & $c$ \\
                 \hline
                 \hline
\textrm{polynomial}       & $4$   & $-0.8076923$ & -- & -- \\
                 & $5$   & $23.636752$  & $-9.(7)$ & -- \\
                 & $6$   & $6.956942$ & $5.627256$ & $-5.135011$        \\
                 \hline
                 \hline
\textrm{trigonometric}    & $4$   & $0.0202$ & --   & --
                 \\
			     & $5$   & $0.785988$ & $-0.356639$& --
                 \\
                 & $6$   & $1.0407$   & $-0.312242$& $-0.0105136$
                 \\
                 \hline
\end{tabular}
\caption{\label{tab1} Values of the coefficients $a$, $b$, and $c$ for the polynomial and trigonometric expansions of the STA solution (see the Supplementary Sec.~I for a complete description).}
\end{table}

We finally assume that the trajectory can be expressed as the sum of some exponential real functions:
\begin{equation}
\label{expeq}
x(t)=a \e^t+b \e^{-t}+c \e^{k t}+d \e^{-k t},
\end{equation}
where $k$ is a free real parameter. This ansatz which is inspired by optimal control theory will become clear in the next section.
From the boundary conditions $x(0)=0, x(T=1)=1$ and $\dot{x}(0)=\dot{x}(T=1)=0$, we obtain:
\begin{equation}
B \vec{h}=\vec{s},
\end{equation}
where
\begin{equation}
B=
\left(
\begin{array}{cccc}
1 & 1 & 1 & 1
\\
\e & \e^{-1} & \e^k & \e^{-k}
\\
1 & -1 & k & -k
\\
\e & -\e^{-1} & k\e^k & -k\e^{-k}
\end{array}
\right),
\end{equation}
and $\vec{h}=(a,b,c,d)^\intercal$, $\vec{s}=(0,1,0,0)^\intercal$. We deduce that $\vec{h}=B^{-1}\vec{s}$. The complete analytical expression of $\vec{h}$ is given in Supplementary Sec.~I.

Table~\ref{tab2} gives the cost functional $\mathcal{C}$ for different STA protocols and for the global minimum solution which is derived in Sec.~\ref{sec4} using optimal control techniques. As could be expected for polynomial and trigonometric expansions, the higher the order $N$ is, the smaller $\mathcal{C}$. We observe nevertheless that these solutions remain quite far from the global minimum solution, 6\% and 12\% respectively for the polynomial and trigonometric solutions at the order 6. The basis of real exponential functions seems to be more suited to the control problem since a small cost (1\% larger than the optimal one) is achieved with the simple expansion proposed in Eq.~\eqref{expeq}.
\begin{table}
\centering
\begin{tabular}{|c|c|c|}
\hline
                 & \textrm{Order} &  \textrm{Cost functional}
                 \\
                 \hline
                 \hline
                 & $N$ 	 & $\mathcal{C}$
                 \\
                 \hline
\textrm{optimal solution} & --    & $1.3130$
\\
\hline
 		         & $3$   & $1.57143$ \\
\textrm{polynomial}       & $4$   & $1.55797$ \\
                 & $5$   & $1.40276$                  \\
                 & $6$   & $1.39986$                  \\
                 \hline
                 & $3$   & $1.70041$                  \\
\textrm{trigonometric}    & $4$   & $1.69843$                  \\
			     & $5$   & $1.48104$                  \\
                 & $6$   & $1.48099$                  \\
                 \hline
\textrm{exponential}      & --    & $1.325271$ \\
\hline
\end{tabular}
\caption{\label{tab2} Cost functional $\mathcal{C}$ for the optimal and the polynomial, trigonometric and exponential STA solutions. The parameter $k$ is set to 100 for the exponential functions. The optimal field is the singular one derived in Sec.~\ref{sec4}.}
\end{table}

\section{The singular geometric optimal solution}\label{sec4}
The control problem we consider belongs to a general class of linear optimal control problems for which powerful mathematical tools, such as the PMP, have been developed~\cite{pont,bryson,bonnardbook}. We solve this problem in this paragraph in a very general setting without any constraint on the control field. In this case, we will show that singular solutions with unbounded fields minimize the cost functional $\mathcal{C}$. Singular control fields have been exhibited in quantum physics in different examples~\cite{lapert:2010,lapert:2013,wu:2012,riviello:2014}.

Before applying optimal control techniques, we need to reformulate the dynamics of the system. To accommodate for the extra boundary conditions on the first time derivative of the speed, we increase the number of state variables over which the optimization is performed by introducing the variable $z$. We will discuss in Sec.~\ref{sec6} a generalization of this strategy to ensure the robustness of optimal solutions to higher derivative orders. We introduce the variables $y=\dot{x}$ and $z=u$. We have:
$$
\ddot{x}+\dot{x}=\dot{u}=v,
$$
where $v(t)$ is the new control field. The dynamical system is now defined on $\mathbb{R}^2$ by:
\begin{equation}
\begin{cases}
\dot{y}=v-y \\
\dot{z}=v,
\end{cases}
\end{equation}
with the boundary conditions $y(0)=y(T)=0$ and $z(0)=0$, $z(T)=1$. The original variables are determined by the relations $x=z-y$ and $u=z$. The cost functional $\mathcal{C}$ to minimize is here given by:
$$
\mathcal{C}=\int_0^T[(z-y)^2+y^2]dt.
$$
The description of the dynamics is by now well suited to the application of the PMP. However, since there is no direct constraint on the amplitude of the control field $v(t)$, the optimal control problem is singular. We assume that the optimal solution is of the form $B-S-B$ where $B$ is a bang pulse and $S$, a singular one. The bang pulse is here a Dirac pulse because in this ideal situation the amplitude of the field is not bounded. The structure of the optimal field will be confirmed by the regularization process introduced in Sec.~\ref{sec5}.

The singular Pontryagin Hamiltonian $H_S$ can be expressed as:
$$
H_S=p_y(v-y)+p_zv-\frac{1}{2}[(z-y)^2+y^2],
$$
where $p_y$ and $p_z$ are the adjoint states associated respectively to $y$ and $z$. We deduce that the Hamilton's equations are:
\begin{equation}
\begin{cases}
\dot{p}_y=-\frac{\partial H_S}{\partial y}=p_y-z+2y \\
\dot{p}_z=-\frac{\partial H_S}{\partial z}=z-y.
\end{cases}
\end{equation}
The singular control fulfills the constraint:
$$
\frac{\partial H_S}{\partial v}=0,
$$
which leads to the definition of the singular set: $p_y+p_z=0$. In the singular case, this relation is satisfied for a non-zero time interval, so that the time derivatives of $p_y+p_z$ are also equal to zero. From $\dot{p}_y+\dot{p}_z=0$, we arrive at $p_y=-y$. From $\dot{p}_y+\dot{y}=0$, we obtain $v=z$. Therefore, we arrive for the singular trajectory $(y_s(t),z_s(t))$ at:
\begin{equation}
\label{eqsing}
\begin{cases}
y_s=Ye^{-t}+Z\sinh (t)\\
z_s=Ze^t.
\end{cases}
\end{equation}
The singular trajectories are labeled by the two parameters $Y$ and $Z$. We observe from Eq.~\eqref{eqsing} that the initial point cannot belong to the singular set, showing the necessity of a bang pulse at time zero. The bang control is a Dirac pulse of amplitude $v_\tau$ and of very small duration $\tau$ such that $v_\tau \tau=\mathcal{A}$ where $\mathcal{A}$ is the area of the pulse, which remains finite in the limit $\tau\to 0$. We denote by $\mathcal{A}_1$ and $\mathcal{A}_2$ the areas of the first and second bangs, respectively. During the bang pulse of very large amplitude, the dynamical system is governed by:
\begin{equation}
\begin{cases}
\dot{y}=v \\
\dot{z}=v,
\end{cases}
\end{equation}
in which the drift term $(y,0)^\intercal$ has been neglected. Since $y(0)=z(0)=0$, we deduce that $y(0^+)=z(0^+)=\mathcal{A}_1$. Note that $t=0^+$ corresponds to the time right after the first bang. Finally, we arrive at the following linear system:
\begin{equation}
\begin{cases}
\mathcal{A}_1(e^{-T}+\sinh(T))+\mathcal{A}_2=0 \\
\mathcal{A}_1e^T+\mathcal{A}_2=1.
\end{cases}
\end{equation}
The solutions can be expressed as:
\begin{equation}
\begin{cases}
\mathcal{A}_1=\frac{1}{\sinh(T)} \\
\mathcal{A}_2=-\frac{1}{\tanh(T)}.
\end{cases}
\end{equation}
We deduce that $z_s(t)-y_s(t)=x_s(t)=\frac{\sinh(t)}{\sinh(T)}$ and $y_s(t)=\frac{\cosh(t)}{\sinh(T)}$. For $T=1$, we have $\mathcal{C}_S=1.3130$, which is given in Tab.~\ref{tab2}. Note that the bang pulses are not taken into account in the computation of the cost functional. This solution is the global minimum of the control problem. We point out that Dirac pulses at time boundaries are required to fulfill the boundary conditions on the first time derivative of the speed. In the absence of such boundary conditions, the optimal solution boils down to the very same singular solution derived above. As a matter of fact, the Dirac pulses originate from a mathematical limit but may not be physically relevant. Optimal control laws of finite amplitude can be derived through a regularization of the cost functional which is described in Sec.~\ref{sec5}.
\section{Regular optimal solutions}\label{sec5}
\subsection{Energy regularization}\label{secrega}
The singular solution derived in Sec.~\ref{sec4} is interesting since it gives the physical limit of the process under study. As shown in Tab.~\ref{tab1} and \ref{tab2}, it also allows us to quantify the efficiency of other solutions of the control problem. However, due to its unbounded character, the corresponding field may pose a problem for its experimental implementation. We propose in this paragraph a regularization of the cost functional~\cite{bryson,liberzon}, denoted $\mathcal{C}_R$, and defined by:
$$
\mathcal{C}_R=\int_0^T[(z-y)^2+y^2+\lambda v^2]dt,
$$
where $\lambda$ is a positive parameter. The singular case is obtained in the limit $\lambda\mapsto 0$. Note that other regularizations such as a bound on the maximum intensity of the field could be also considered. The energy regularization has the advantage to generate a smooth control field, its maximum amplitude being connected to the parameter $\lambda$.

In this case, the regular Pontryagin Hamiltonian $H_R$ can be written as:
$$
H_R=H_S-\frac{\lambda}{2}v^2.
$$
We have the same Hamilton's equations as in the singular situation but the regular control field is given by:
$$
v_R=\frac{p_y+p_z}{\lambda}.
$$
We deduce the following relations:
\begin{equation}\label{eqreg}
\begin{cases}
\lambda \dot{y}=p_y+p_z-\lambda y \\
\lambda \dot{z}=p_y+p_z \\
\dot{p}_y=p_y-z+2y \\
\dot{p}_z=z-y .
\end{cases}
\end{equation}
The final step consists in integrating these equations and finding the initial adjoint states $p_y(0)$ and $p_z(0)$ so that the final boundary conditions are satisfied: $y(T)=0$ and $z(T)=1$. A general theory called linear quadratic optimal control was developed to solve this kind of optimal control problems~\cite{bryson,bressan,liberzon}. We recall in Sec.~\ref{secregb} the basic features of this theory.
\subsection{Linear quadratic optimal control}\label{secregb}
We show in this paragraph that the optimal solution of some specific linear quadratic optimal control problems can be expressed in terms of real or complex exponential terms. Here, we briefly outline the different steps of this approach. The interested reader can find rigorous mathematical derivations in some standard textbooks~\cite{bryson,bressan,liberzon}.

The control problem is defined as follows. The state of the system is a vector $x\in\mathbb{R}^n$ whose dynamics are governed by the following differential equation:
\begin{equation}
\dot{x}=Ax+Bu,
\label{ABequation}
\end{equation}
where the control field $u(t)\in\mathbb{R}^m$, $A\in M_n(\mathbb{R})$ and $B\in M_{nm}(\mathbb{R})$ are two constant matrices. Starting from the state $x(0)=x_0$, the goal is to reach the state $x(T)=x_f$ at time $T$, while minimizing the cost functional $\mathcal{C}$ defined by:
\begin{equation}
\mathcal{C}=\int_0^T[x^\intercal W x+u^\intercal U u]dt,
\end{equation}
where $W$ and $U$ are two constant symmetric matrices which are respectively positive and positive definite. The Pontryagin Hamiltonian $H_P$ of the system can be written as:
\begin{equation}
H_P=p^\intercal Ax+p^\intercal Bu-\frac{1}{2}(x^\intercal W x+u^\intercal Uu),
\end{equation}
where $p\in\mathbb{R}^n$ is the adjoint state. The differential equation governing the dynamics of $p$ is:
\begin{equation}
\dot{p}=-A^\intercal p+Wx.
\end{equation}
The optimality condition $\frac{\partial H_P}{\partial u}=0$ leads to:
\begin{equation}
u=U^{-1}B^\intercal p.
\end{equation}
We finally obtain the optimal equations for the vectors $x$ and $p$:
\begin{equation}\label{eqesystem}
\begin{pmatrix}
\dot{x} \\
\dot{p}
\end{pmatrix}
=
\begin{pmatrix}
A & BU^{-1}B^\intercal \\
W & -A^\intercal
\end{pmatrix}
\begin{pmatrix}
x \\ p
\end{pmatrix}
\end{equation}
The solutions of this first order linear differential system can be expressed as the sum of exponential functions. Using Eq.~\eqref{eqesystem}, the last step of the method consists in finding the initial adjoint state $p(0)$ such that the corresponding trajectory reaches the target $x_f$ at time $t$. This can be made by deriving the solution of Eq.~\eqref{eqesystem}. We emphasize that the previous reasoning can be readily generalized to time-dependent matrices $A$ and $B$ in Eq.~(\ref{ABequation})~\cite{bryson}.

In the example under study in Sec.~\ref{secrega}, it is straightforward to show that:
\begin{equation}
A=\begin{pmatrix}
-1 & 0 \\
0 & 0
\end{pmatrix},~B=\begin{pmatrix} 1 \\ 1\end{pmatrix},~W=\begin{pmatrix}
2 & -1 \\
-1 & 1
\end{pmatrix}, U=\lambda.
\end{equation}
Note that $W$ is positive since its two eigenvalues are real positive numbers. For this example, Eq.~\eqref{eqreg} are equivalent to Eq.~\eqref{eqesystem}. The eigenvalues of the matrix defined in Eq.~\eqref{eqesystem} are 1, -1, $1/\sqrt{\lambda}$ and $-1/\sqrt{\lambda}$, so the time evolution of the states and adjoint states can be expressed as the sum of real exponential functions as in Eq.~\eqref{expeq}, with $k=1/\sqrt{\lambda}$ (see Supplementary Sec. III for details). This point explains the efficiency of the corresponding STA process which corresponds to a regularized optimal solution. It also shows how OCT can help selecting the best basis of functions to expand the control field or the trajectory of the system.
\section{Comparison between STA and optimal control protocols}
\begin{figure}[htp]
\centering
\includegraphics[width=0.75\textwidth]{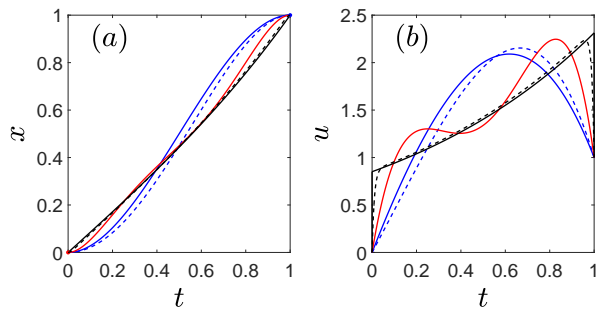}
\caption{(Color online) Comparison between the STA and the optimal solutions. Panels (a) and (b) represent respectively the time evolution of the trajectory $x$ and of the control field $u$. The solid and dashed black lines correspond respectively to the singular and regular solutions. The parameter $\lambda$ is set to $10^{-4}$. The STA polynomial solutions are also plotted in blue (solid and dashed) and red for the orders $n=3$, 4 and 5. The solution at the order 6 is very similar to the one at the order 5. Dimensionless units are used.}
\label{fig1}
\end{figure}

\begin{figure}[htp]
\centering
\includegraphics[width=0.75\textwidth]{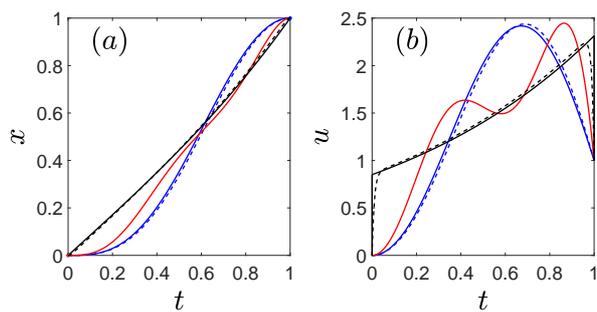}
\caption{(Color online) Same as Fig.~\ref{fig1} but for the trigonometric STA solution. Dimensionless units are used.}
\label{fig2}
\end{figure}

We present in this section different numerical results to illustrate the comparison between STA and OCT solutions and between singular and regular optimal fields. Figures \ref{fig1} and \ref{fig2} display the trajectory and the control field corresponding to the STA and to the singular processes. In the two cases, we observe that the STA solutions present strong oscillations and are not able to reproduce optimal dynamics, even for high order polynomial or trigonometric expansion. This point is a clear indication that polynomial and trigonometric functions are not a good basis to expand the control solution. This conclusion confirms the results of the cost functional presented in Tab.~\ref{tab1} and~\ref{tab2}.

\begin{figure}[h!]
\centering
\includegraphics[width=0.75\textwidth]{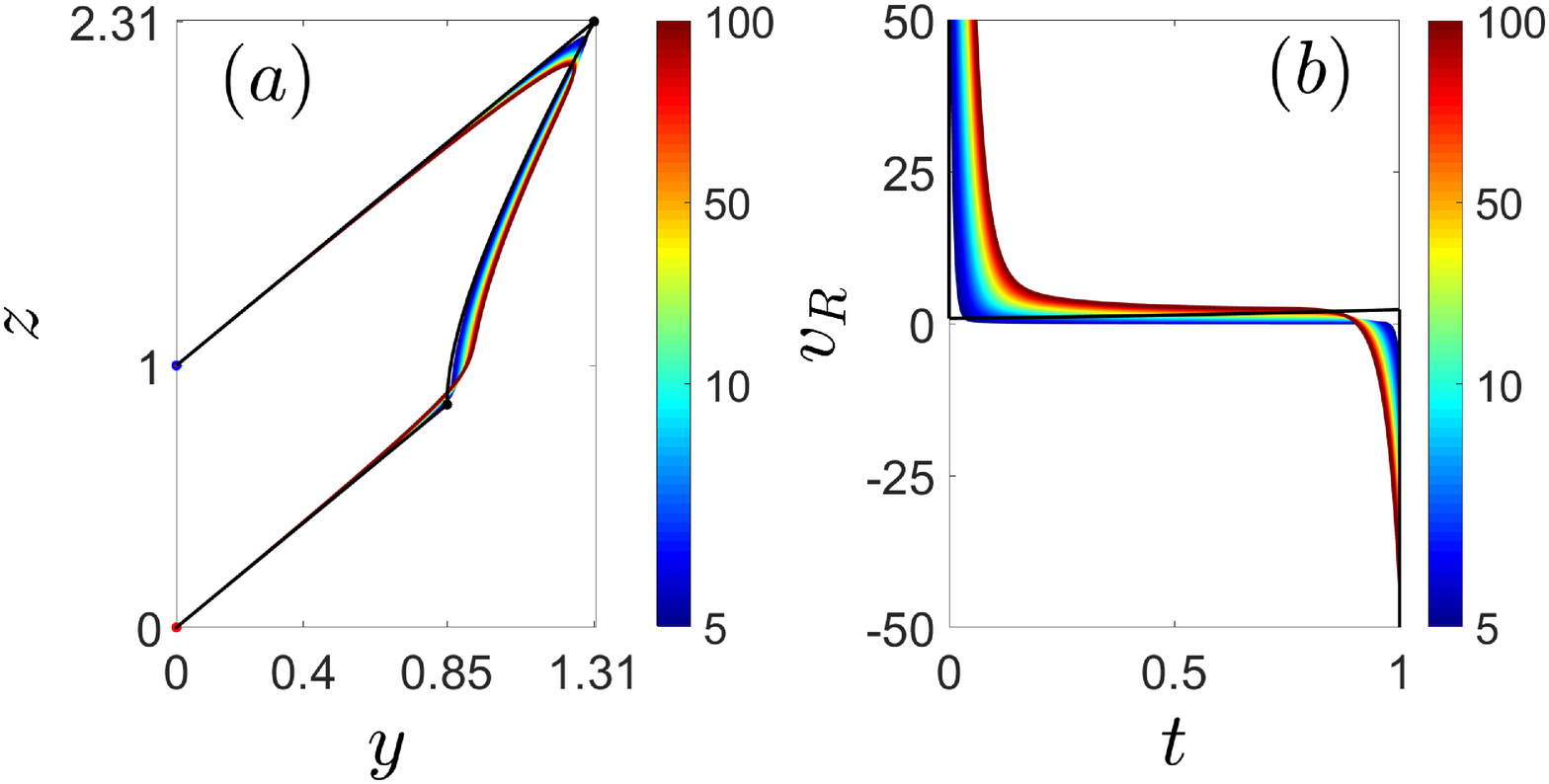}
\caption{(Color online) Panel (a): Plot of the optimal regular and singular trajectories in the $(y,z)$- plane. The dots indicate the initial and final points of the trajectories corresponding to bang pulses. Panel (b): Plot of the time evolution of the regular control field $v_R$. The black solid line corresponds to the singular solution, while the regular solutions are displayed in color. The color bar indicates the value of $\lambda\times 10^5$ for each regular process. Dimensionless units are used.}
\label{fig3}
\end{figure}

\begin{figure}[htp]
\centering
\includegraphics[width=0.75\textwidth]{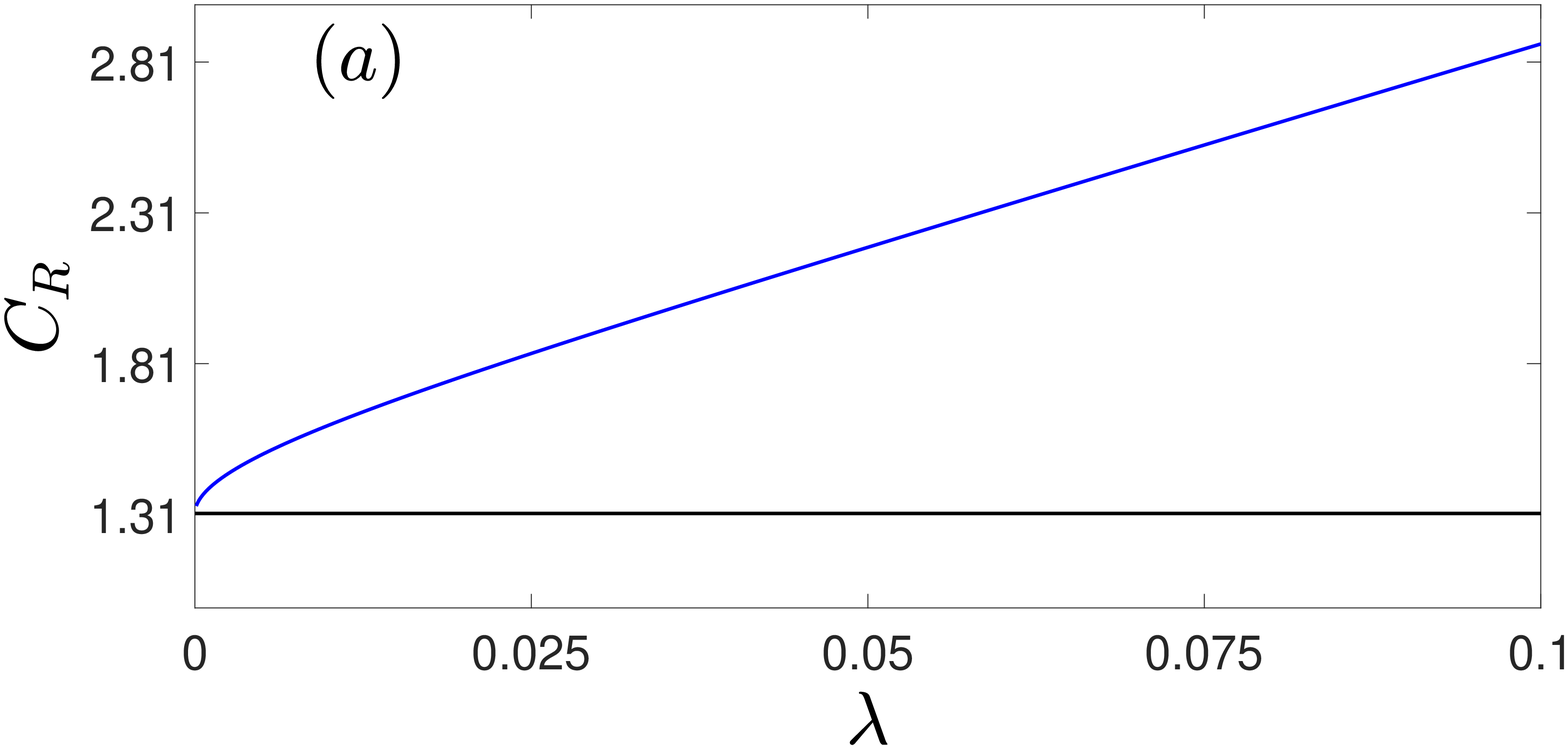}
\includegraphics[width=0.75\textwidth]{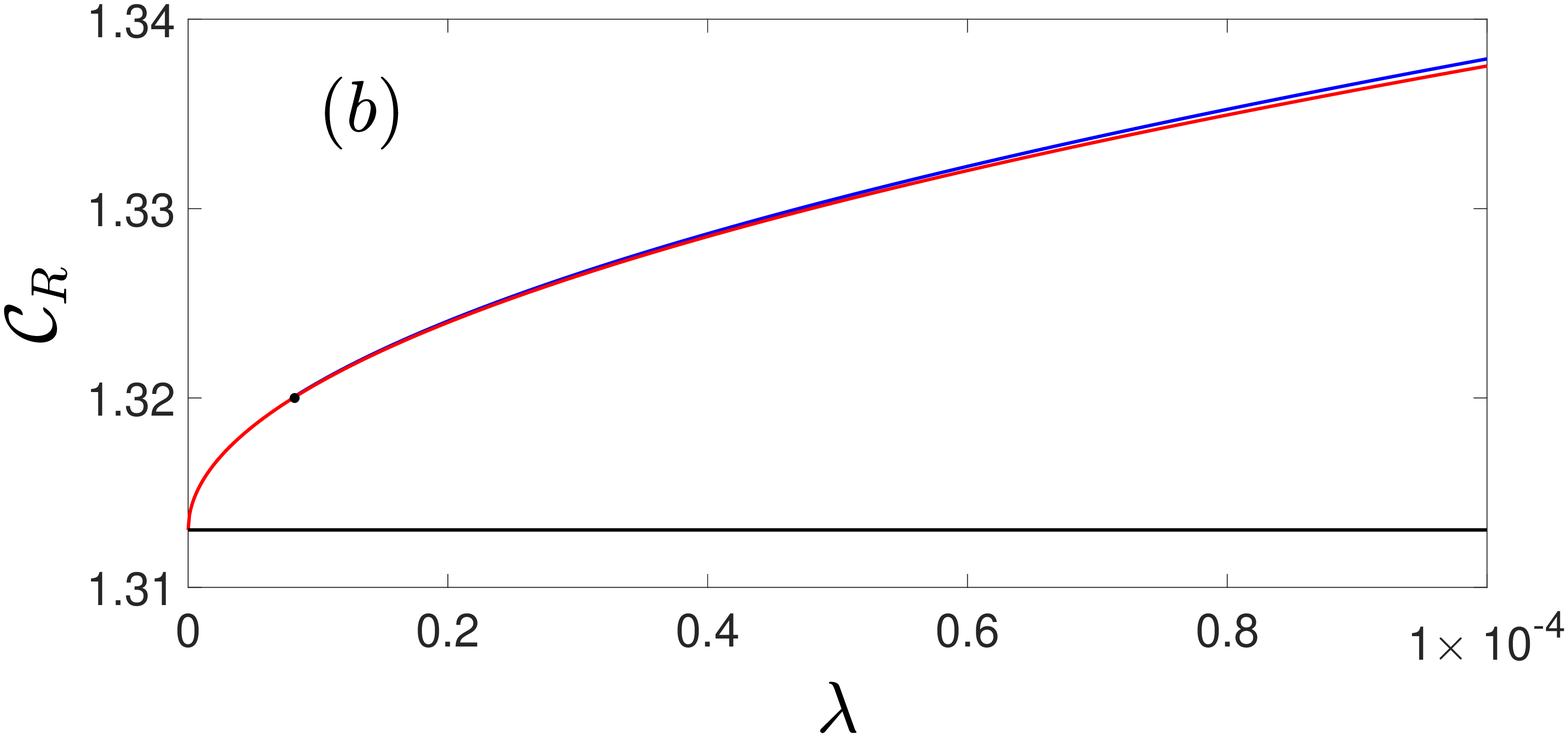}
\caption{(Color online) Panel (a): Plot of the regularized cost functional $C_R$ as a function of the parameter $\lambda$ (solid blue line). The solid horizontal black line indicates the value of the minimum cost functional for the singular solution. Panel (b) is a zoom of panel (a) for $0\leq \lambda \leq 10^{-4}$. The red curve is a square root fit around $\lambda=0$ of $C_R$. The black dot indicates the minimum value of $\lambda$ for which $C_R$ has been computed numerically. Dimensionless units are used.}
\label{fig4}
\end{figure}

Figures~\ref{fig3} and \ref{fig4} illustrate the convergence of the regular optimal solution towards the singular one when the parameter $\lambda$ goes to 0. Figure~\ref{fig3} shows that the regular trajectory and control field converge smoothly towards the singular ones. This numerical observation is an evidence strengthening the conjecture made in Sec.~\ref{sec4} about the structure of the optimal solution. The evolution of the cost functional $C_R$ as a function of $\lambda$ is displayed in Fig.~\ref{fig4}. Numerical simulations reveal a rapid convergence to the minimum value of 1.3130. As shown in Fig.~\ref{fig4}b, $C_R$ exhibits a square root behavior in the limit $\lambda\to 0$. Due to the complexity of the computations, we were not able to prove this statement analytically. From a numerical point of view, the regular solution can be derived up to $\lambda=2\times 10^{-6}$, divergences occur for lower values of $\lambda$.

\section{Generalization to higher derivative orders}\label{sec6}
This section is aimed at generalizing the control problem to constraints of any order on the initial and final time derivatives of the trajectory. The control protocol is defined  by the differential equation $\dot{x}+x=u$ with the boundary conditions $x(0)=\dot{x}(0)=\cdots =x^{(n)}(0)=0$ and $x(T)=1$, $\dot{x}(T)=\cdots =x^{(n)}(T)=0$. These conditions at the boundary of the control time interval lead to a robustness of the protocol against slight variations at initial and final times. This property is required to ensure a reliable experimental implementation of control processes. As in Sec.~\ref{sec4}, we enlarge the space of variables for the optimization procedure by introducing extra coordinates $x_1=\dot{x}$, $x_2=\ddot{x}$, $\cdots$, $x_n=x^{(n)}$ and $z_0=u$, $z_1=\dot{u}$, $\cdots$, $z_{n-1}=u^{(n-1)}$ to replace the boundary conditions on the derivatives by conditions on the state of the system. The new control term is $v=u^{(n)}$. The dimension of the state of the system is $n+1$, $(x_n,z_{n-1},\cdots,z_1,z_0)$. The differential system to control can be expressed as:
\begin{equation}
\begin{cases}
\dot{x}_n+x_n=v \\
\dot{z}_{n-1}=v \\
\dot{z}_{n-2}=z_{n-1}\\
\cdots \\
\dot{z}_k=z_{k+1} \\
\cdots \\
\dot{z}_1=z_2 \\
\dot{z}_0=z_1,
\end{cases}
\end{equation}
with the boundary conditions $x_n(0)=x_n(T)=0$, $z_{n-1}(0)=z_{n-1}(T)=0$, $\cdots$, $z_1(0)=z_1(T)=0$ and $z_0(0)=0$, $z_0(T)=1$. The cost functional in the regular case can be written as:
$$
\mathcal{C}_R=\int_0^T[(z_0-x_1)^2+x_1^2+\lambda v^2]dt,
$$
with
$$
x_1=z_1-z_2+z_3+\cdots +(-1)^nz_{n-1}-(-1)^nx_n.
$$
The singular limit is obtained for $\lambda =0$.

\begin{figure}[htp]
\centering
\includegraphics[width=0.75\textwidth]{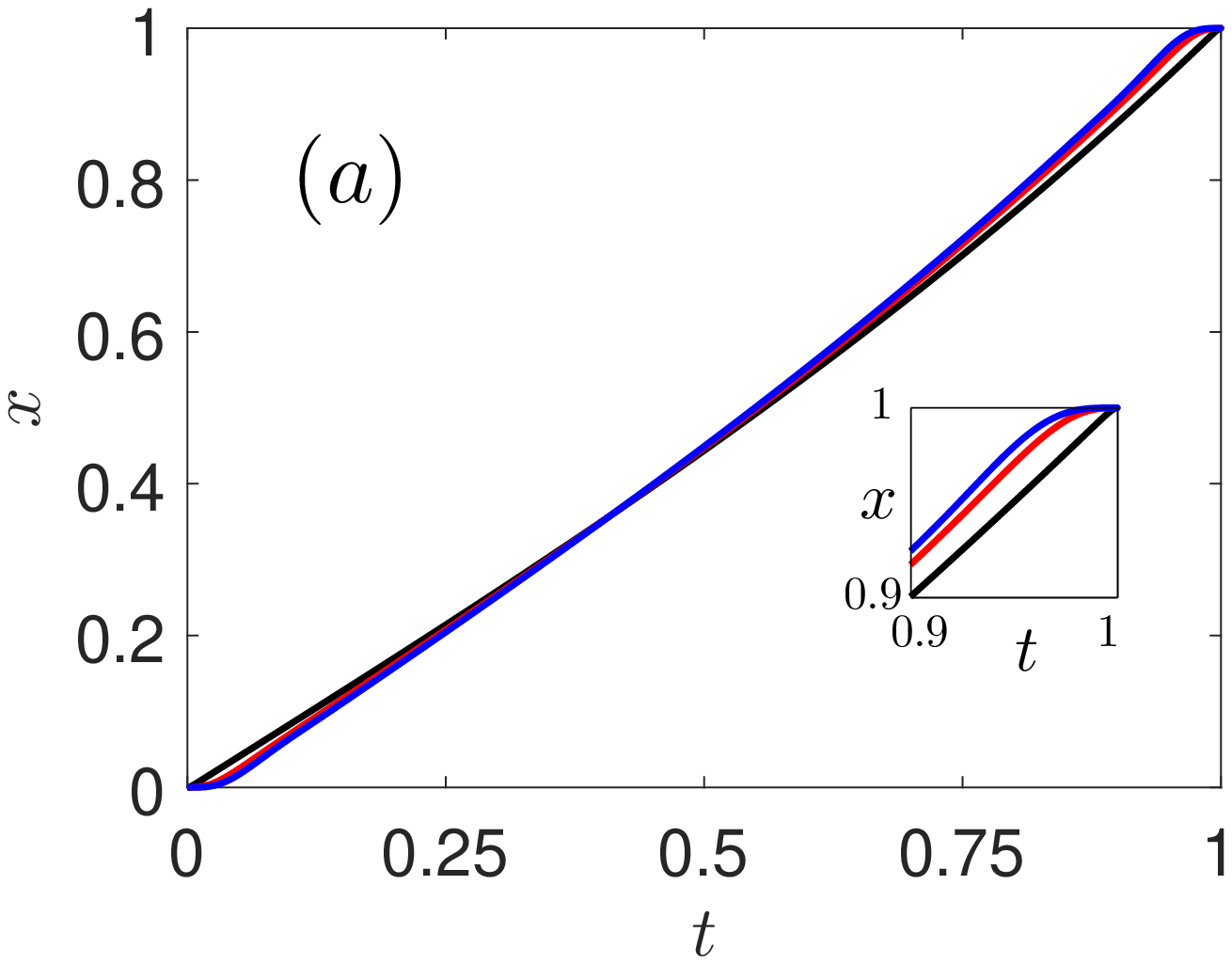}
\includegraphics[width=0.75\textwidth]{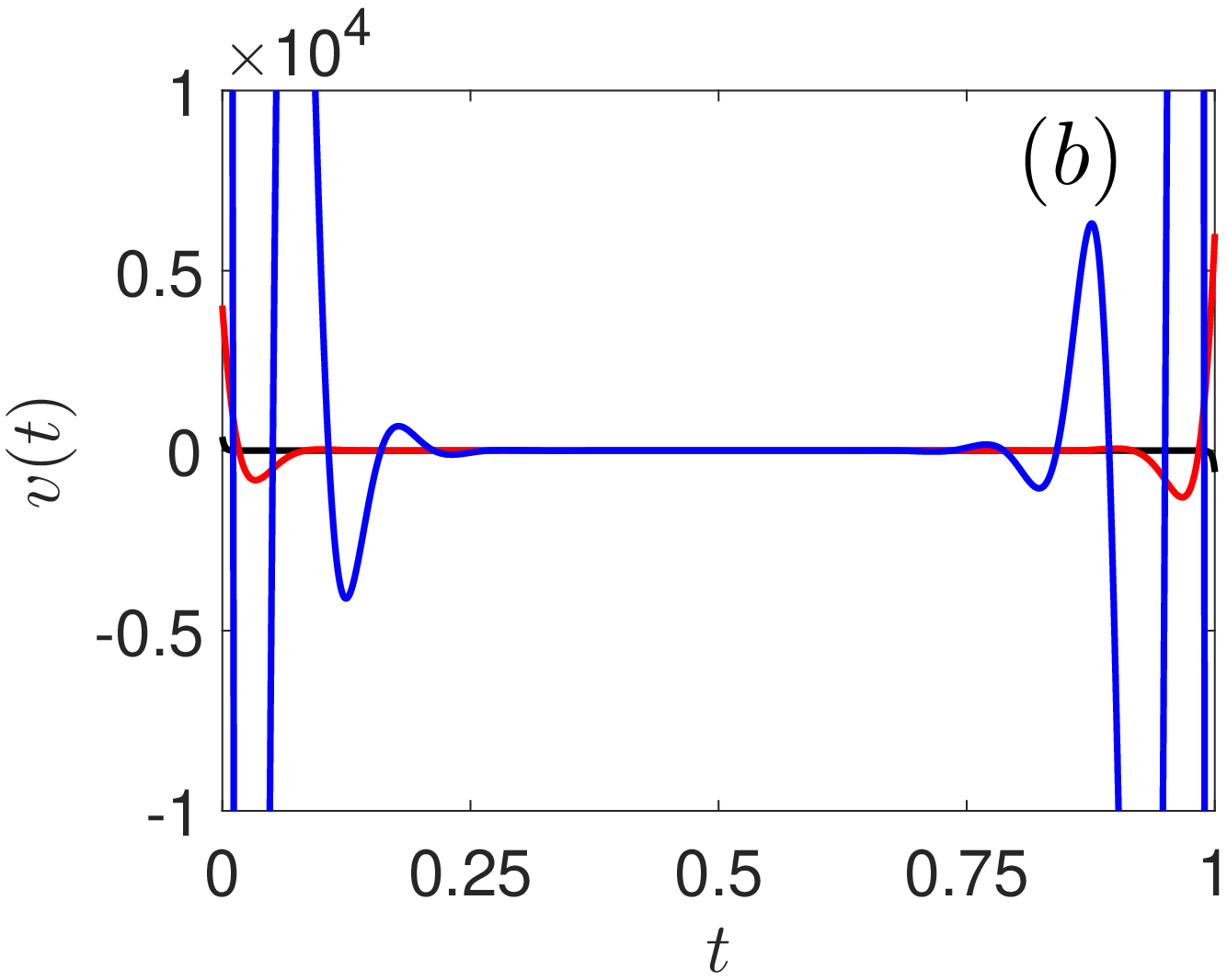}
\caption{(Color online) Plot of the time evolution of the optimal trajectory (a) and of the regular control field (b) for $n=1$ (black), $n=2$ (red) and $n=3$ (blue). The parameter $\lambda$ is set respectively to $10^{-5}$, $5\times 10^{-7}$ and $5\times 10^{-9}$ for $n=1$, 2 and 3. In panel (a), the small insert is a zoom of the trajectories near $t=T$. Dimensionless units are used.}
\label{fig5}
\end{figure}

The regular control can be derived by using the material of Sec.~\ref{secregb}. For each order, it can be verified that the matrices $W$ and $U$ satisfy the constraints of the method. Different regular solutions for $n=1$ to 3 are plotted in Fig.~\ref{fig5}. The behavior of the different trajectories at time 0 and $T=1$ is clearly seen in Fig.~\ref{fig5}. The case $n=1$ is the only case where the matrix defined in Eq.~\eqref{eqesystem} has real eigenvalues. For $n>1$, some of the eigenvalues are complex numbers and the trajectory is expressed as a linear combination of real exponential and trigonometric functions. The arguments of the exponential terms can be determined explicitly from Eq.~\eqref{eqesystem}. The optimal trajectory and control field have a relatively complicated evolution at the beginning and at the end of the control process, due to the different boundary conditions to fulfill. However, we observe a similar behavior at intermediate times for the three orders. As in the case $n=1$, this corresponds to the singular set of the optimal control problem, which is the same at any order.
This statement can be proved for any order $n$. For that purpose, we introduce the singular Pontryagin Hamiltonian which can be expressed as:
\begin{eqnarray*}
& & H_S=p_n(v-x_n)+p_{n-1}v+p_{n-2}z_{n-1}+\cdots +p_0z_1\\
& & -\frac{(z_0-x_1)^2+x_1^2}{2},
\end{eqnarray*}
where the adjoint states $p_n$, $p_{n-1}$, ..., $p_1$ and $p_0$ are associated respectively to $x_n$, $z_{n-1}$, ..., $z_1$ and $z_0$. Using $H_S$, we also deduce the differential equations governing the dynamics of the $p_i$s:
\begin{equation}
\begin{cases}
\dot{p}_n=p_n+(-1)^{n+1}(2x_1-z_0) \\
\dot{p}_{n-1}=-p_{n-2}+(-1)^n(2x_1-z_0) \\
\dot{p}_{n-2}=-p_{n-3}+(-1)^{n-1}(2x_1-z_0)\\
\cdots \\
\dot{p}_k=-p_{k-1}+(-1)^{k+1}(2x_1-z_0) \\
\cdots \\
\dot{p}_1=-p_0+2x_1-z_0 \\
\dot{p}_0=z_0-x_1.
\end{cases}
\end{equation}
The singular set is characterized by $\frac{\partial H_S}{\partial v}=0$, which leads to $p_n+p_{n-1}=0$. If the trajectory belongs to this set in a non zero time interval then the time derivatives of $p_n+p_{n-1}$ are also equal to 0. We then obtain a series of constraints given by:
\begin{equation}
\begin{cases}
p_n-p_{n-2}=0 \\
p_n+p_{n-3}=0 \\
\cdots \\
p_n+(-1)^{n-1}p_0=0\\
p_n+(-1)^{n+1}x_1=0 \\
z_1=z_0
\end{cases}
\end{equation}
where each relation is obtained by deriving with respect to time the preceding one. From $z_1=z_0$, we deduce that $z_0=z_1=z_2=\cdots =z_{n-1}=v$ on the singular set. The singular trajectories can then be written as:
\begin{equation}
\begin{cases}
x_n(t)=Ye^{-t}+Z\sinh(t) \\
z_0(t)=\cdots =z_{n-1}(t)=v(t)=Ze^t
\end{cases}
\end{equation}
where the two constants $Y$ and $Z$ characterize the singular solutions. Note that the singular control field is the same for any value of $n$.

\section{Conclusion and prospective views}\label{sec7}

Optimal Control Theory protocols are built on about a global constraint, the minimization of a cost functional, while STA techniques are primarily built to accommodate for local constraints, in particular at time interval boundaries. The thorough analysis of the control of a linear system with boundary conditions at initial and final times has enabled us to establish how the two formalisms can mutually benefit from each other.

In the simple STA approach used in this manuscript, the trajectory is expanded over a basis of functions and depends on a finite number of parameters. Such parameters are determined to fulfill the boundary conditions and to minimize the cost functional. We have shown that polynomial and trigonometric functions are not well suited to the control problem under study, while a basis of real exponential functions inspired by regularized optimal solutions gives very efficient results.

Interestingly, STA protocols can be made robust against uncertainties on the initial and final times by the cancellation of successive time derivatives at time interval boundaries. In contrast, usual OCT solutions are perfectly suited for a well defined given time interval but fail down in presence of time interval uncertainties. Inspired by STA local constraints, we have improved the robustness of OCT solutions by enlarging the parameter space of the optimal control problem so to accommodate for local constraints at time interval boundaries.

This study can also be viewed as a test case showing the efficiency and the flexibility of linear quadratic optimal control theory. This approach, well-known in mathematics and in engineering, has been little applied in physics. Such methods can be used when the system is governed by a linear differential equation such as, e.g., in ion cyclotron resonance~\cite{bodenhausen:2016,delsuc:2013,delsuc:2016}, but also in a non-linear framework in a range where the linear approximation is valid~\cite{li:2017}.

\noindent\textbf{ACKNOWLEDGMENT}\\
D. Sugny acknowledges support from the PICS program, from the ANR-DFG research program COQS (ANR-15-CE30-0023-01) and from the QUACO project (ANR 17-CE40-0007-01). This project has received funding from the European Union's Horizon 2020 research and innovation programme under the Marie-Sklodowska-Curie grant agreement No 765267 (QUSCO). This work has been financially supported by the Agence Nationale de la Recherche research funding Grant No. ANR-18-CE30-0013.

\newpage

\appendix

\textbf{\Large Supplementary material of the paper:\\
A comparison between optimal control and shortcut to adiabaticity protocols in a linear control system}\\ \\

This supplementary material gives technical details about the derivation of STA and regular optimal solutions in Sec. A and B, respectively.\\ \\

\section{STA protocols}
We start from the differential equation:
\begin{equation}
\dot{x}+x=u,
\end{equation}
where $x(t)$ and its derivative satisfy the boundary conditions $x(0)=0$, $x(T)=1$ and $\dot{x}(0)=\dot{x}(T)=0$. We consider different function families in which the trajectory of the system $x(t)$ is expanded.
\subsection{Polynomial function family}
We assume that
\begin{equation}
\label{xp}
x(t)=\sum\limits_{k=2}^N a_k t^k
\end{equation}
From $x(T)=1$ and $\dot{x}(T)=0$, we have:
\begin{equation}
\label{xp1}
\begin{cases}
\sum\limits_{k=2}^N a_k T^k=1
\\
\sum\limits_{k=2}^N k a_k T^{k-1}=0
\end{cases}
\end{equation}
So, there are 2 equations and $N-1$ unknowns, i.e. $N-3$ free parameters.\\
\noindent $N=3$:\\
Equation~\eqref{xp} becomes:
\begin{equation}
x_3(t)=a_2 t^2+a_3 t^3
\end{equation}
Equation~\eqref{xp1} leads to the system:
\begin{equation}
\begin{cases}
T^2 (a_2+a_3 T)=1
\\
T (2 a_2+3 a_3 T)=0
\end{cases}
\end{equation}
which gives
\begin{equation}
a_2=\frac{3}{T^2},\quad a_3=-\frac{2}{T^3},
\end{equation}
and we obtain:
\begin{equation*}
x_3=3\left(\frac{t}{T}\right)^2-2\left(\frac{t}{T}\right)^3
\end{equation*}
The control field takes the form:
\begin{gather*}
u_3=\frac{6}{T}\left(\frac{t}{T}\right)+\frac{3}{T} (T-2)\left(\frac{t}{T}\right)^2-2\left(\frac{t}{T}\right)^3
\end{gather*}
The cost functional has been numerically evaluated as $\mathcal{C}_3\approx 1.57143$, when $T=1$.\\
\noindent $N=4:$\\
We have:
\begin{equation}
x_4=a_2 t^2+a_3 t^3+a_4 t^4,
\end{equation}
and using the boundary conditions, we deduce that:
\begin{equation*}
a_2=\frac{3+T^4 a}{T^2},\quad a_3=-2\frac{1+T^4 a}{T^3},\quad a_4=a,
\end{equation*}
where $a$ is a free parameter. We arrive at:
\begin{equation*}
x_4(t)=(3+T^4 a) \left(\frac{t}{T}\right)^2-2 (1+T^4 a) \left(\frac{t}{T}\right)^3+T^4 a \left(\frac{t}{T}\right)^4
\end{equation*}
The minimization in this case can be done analytically. Plugging the solution into the cost functional yields:
\begin{equation*}
\mathcal{C}_4=2\left[\left(\frac{1}{1260}+\frac{1}{105}\right) a^2+\frac{1}{60} a+\frac{11}{14}\right],
\end{equation*}
and hence
\begin{equation*}
\begin{cases}
a_{min}=-\frac{21}{26}\approx - 0.8076923
\\
\mathcal{C}\approx 1.55797.
\end{cases}
\end{equation*}
\noindent $N=5:$\\
Starting from:
\begin{equation*}
x_5=a_2 t^2+a_3 t^3+a_4 t^4+a_5 t^5,
\end{equation*}
we obtain with the boundary conditions:
\begin{equation*}
a_2=(3+a+2 b),\quad a_3=-(2+2 a+3 b),\quad a_4=a,\quad a_5=b,
\end{equation*}
where $a$ and $b$ are free parameters. We finally arrive for $T=1$ at:
\begin{equation*}
\begin{cases}
x_5=(3+a+2 b) t^2-(2+2 a+3 b) t^3+a t^4+b t^5
\\
u_5=2 (3+a+2 b) t-(3+5 a+7 b) t^2+(-2+2 a-3 b) t^3+(a+5 b) t^4+b t^5
\end{cases}
\end{equation*}
A numerical minimization leads to:
\begin{equation*}
a\approx 23.636752,~ b\approx -9.(7),~ \mathcal{C}_5\approx 1.40276
\end{equation*}
\noindent $N=6:$\\
We have:
\begin{equation*}
x_6=a_2 t^2+a_3 t^3+a_4 t^4+a_5 t^5+a_6 t^6,
\end{equation*}
and from the boundary conditions, we obtain:
\begin{equation*}
a_2=(3+a+2 b+3 c),\quad a_3=-(2+2 a+3 b+4 c),\quad a_4=a,\quad a_5=b,\quad a_6=c
\end{equation*}
where $a$, $b$ and $c$ are free parameters. We get:
\begin{equation*}
x_6=(3+a+2 b+3 c) t^2-(2+2 a+3 b+4 c) t^3+a t^4+b t^5+c t^6
\end{equation*}
and
\begin{align}
\nonumber u_6^\P(t)&=2 (3+a+2 b+3 c) t-(3+5 a+7 b+9 c) t^2+(-2+2 a-3 b-4 c) t^3+
\\
\nonumber &+(a+5 b) t^4+(b+6 c) t^5+c t^6
\end{align}
The corresponding numerical values for the free parameters and for the cost functional are:
\begin{equation*}
a\approx 6.956942,~ b\approx 5.627256,~ c\approx -5.135011,~ \mathcal{C}_6\approx 1.39986.
\end{equation*}
\subsection{Trigonometric function family}
We assume that the trajectory can be expressed as:
\begin{equation}
\label{xtg}
x(t)=\sum\limits_{k=1}^N a_k \sin\left(\frac{k\pi}{2} t\right)
\end{equation}
$x(0)=0$ is automatically satisfied and the condition $\dot{x}(0)=0$ leads to:
\begin{eqnarray}
\sum\limits_{k=1}^N k a_k=0
\end{eqnarray}
From $x(1)=1$ and $\dot{x}(1)=0$, we deduce that:
\begin{equation}
\label{xtg1}
\begin{cases}
\sum\limits_{k=1}^N a_k \sin\left(\frac{k\pi}{2}\right)=1
\\
\sum\limits_{k=1}^N k a_k \cos\left(\frac{k\pi}{2}\right)=0
\end{cases}
\end{equation}
We have 3 equations to fulfill and $N$ unknowns, so $N-3$ free parameters.\\
\noindent $N=3$:\\
Equation~\eqref{xtg} leads to:
\begin{equation}
x_3(t)=a_1 \sin\left(\frac{\pi}{2} t\right)+a_2 \sin\left(\pi t\right)+a_3 \sin\left(\frac{3 \pi}{2} t\right)
\end{equation}
From the boundary conditions, we obtain:
\begin{equation}
a_1=\frac{3}{4},\quad a_2=0,\quad a_3=-\frac{1}{4}.
\end{equation}
We finally arrive at:
\begin{equation*}
\begin{cases}
x_3(t)=\frac{3}{4} \sin\left(\frac{\pi}{2} t\right)-\frac{1}{4}\sin\left(\frac{3\pi}{2} t\right)
\\
u_3(t)=\frac{3}{4}\left(\frac{\pi}{2}\cos\left(\frac{\pi}{2} t\right)+\sin\left(\frac{\pi}{2} t\right)\right)-\frac{1}{4}\left(\frac{3\pi}{2}\cos\left(\frac{3\pi}{2} t\right)+\sin\left(\frac{3\pi}{2} t\right)\right)
\end{cases}
\end{equation*}
The corresponding cost functional is:
\begin{equation}
\mathcal{C}_3\approx 1.70041
\end{equation}
\noindent $N=4$:\\
Starting from:
\begin{equation}
x_4(t)=a_1 \sin\left(\frac{\pi}{2} t\right)+a_2 \sin\left(\pi t\right)+a_3 \sin\left(\frac{3 \pi}{2} t\right)+a_4 \sin\left(2\pi t\right)
\end{equation}
we obtain:
\begin{equation}
a_1=\left(\frac{3}{4}-2 a\right),\quad a_2=2 a,\quad a_3=-\left(\frac{1}{4}+2 a\right),\quad a_4=a,
\end{equation}
where $a$ is a free parameter. Thus we finally arrive at:
\begin{equation*}
x_4(t)=\left(\frac{3}{4}-2 a\right)\sin\left(\frac{\pi}{2} t\right)+2 a \sin\left(\pi t\right)-\left(\frac{1}{4}+2 a\right)\sin\left(\frac{3\pi}{2} t\right)+a \sin\left(2\pi t\right)
\end{equation*}
\begin{align}
\nonumber &u_4(t)=\left(\frac{3}{4}-2 a\right)\left(\frac{\pi}{2}\cos\left(\frac{\pi}{2} t\right)+\sin\left(\frac{\pi}{2} t\right)\right)+2 a\left(\pi\cos\left(\pi t\right)+\sin\left(\pi t\right)\right)-
\\
\nonumber &-\left(\frac{1}{4}+2 a\right)\left(\frac{3\pi}{2}\cos\left(\frac{3 \pi}{2} t\right)+\sin\left(\frac{3\pi}{2}t\right)\right)+a\left(2\pi\cos\left(2\pi t\right)+\sin\left(2\pi t\right)\right)
\end{align}
The corresponding values for the free parameter and for the cost functional are:
\begin{equation}
a\approx 0.0202,~ \mathcal{C}_4\approx 1.69843.
\end{equation}
\noindent $N=5$:\\
Equation~\eqref{xtg} leads to:
\begin{equation}
x_5(t)=a_1 \sin\left(\frac{\pi}{2} t\right)+a_2 \sin\left(\pi t\right)+a_3 \sin\left(\frac{3 \pi}{2} t\right)+a_4 \sin\left(2\pi t\right)+a_5 \sin\left(\frac{5\pi}{2} t\right)
\end{equation}
and we obtain:
\begin{equation}
a_1=\left(\frac{3}{4}-2 a\right),~ a_2=2 (a-b),~ a_3=-\left(\frac{1}{4}+2 a-b\right),~ a_4=a-b,~ a_5=b
\end{equation}
where $a$ and $b$ are free parameters. We get:
\begin{align}
&\nonumber x_5(t)=\left(\frac{3}{4}-2 a\right)\sin\left(\frac{\pi}{2} t\right)+2 (a-b) \sin\left(\pi t\right)-
\\
\nonumber &-\left(\frac{1}{4}+2 a-b\right)\sin\left(\frac{3\pi}{2} t\right)+(a-b) \sin\left(2\pi t\right)+b\sin\left(\frac{5\pi}{2} t\right)
\end{align}
\begin{align}
\nonumber &u_5(t)=\left(\frac{3}{4}-2 a\right)\left(\frac{\pi}{2}\cos\left(\frac{\pi}{2} t\right)+\sin\left(\frac{\pi}{2} t\right)\right)+2 (a-b) \left(\pi\cos\left(\pi t\right)+\sin\left(\pi t\right)\right)-
\\
\nonumber &-\left(\frac{1}{4}+2 a-b\right)\left(\frac{3\pi}{2}\cos\left(\frac{3\pi}{2} t\right)+\sin\left(\frac{3\pi}{2} t\right)\right)+(a-b)\left(2\pi\cos\left(2\pi t\right)+\sin\left(2\pi t\right)\right)+
\\
\nonumber &+b\left(\frac{5\pi}{2}\cos\left(\frac{5\pi}{2} t\right) +\sin\left(\frac{5\pi}{2} t\right)\right)
\end{align}
The corresponding values for the free parameters and for the cost functional are:
\begin{equation}
a\approx 0.785988,~ b\approx -0.356639,~ \mathcal{C}_5\approx 1.48104.
\end{equation}
\noindent $N=6$:\\
The trajectory can be expressed as:
\begin{equation}
x_6(t)=a_1 \sin\left(\frac{\pi}{2} t\right)+a_2 \sin\left(\pi t\right)+a_3 \sin\left(\frac{3 \pi}{2} t\right)+a_4 \sin\left(2\pi t\right)+a_5 \sin\left(\frac{5\pi}{2} t\right)+a_6 \sin\left(3\pi t\right)
\end{equation}
and we obtain:
\begin{equation}
a_1=\left(\frac{3}{4}-2 a-2 b\right),~ a_2=2 a-3 c,~ a_3=-\left(\frac{1}{4}+2 a+b\right),~ a_4=a,\quad a_5=b,~ a_6=c
\end{equation}
where $a$, $b$ and $c$ are free parameters. We finally deduce that:
\begin{align}
&\nonumber x_6(t)=\left(\frac{3}{4}-2 a-2b\right)\sin\left(\frac{\pi}{2} t\right)+(2 a-3 c) \sin\left(\pi t\right)-
\\
\nonumber &-\left(\frac{1}{4}+2 a+b\right)\sin\left(\frac{3\pi}{2} t\right)+a \sin\left(2\pi t\right)+b\sin\left(\frac{5\pi}{2} t\right)+c \sin\left(3\pi t\right)
\end{align}
\begin{align}
\nonumber &u_6(t)=\left(\frac{3}{4}-2 a-2b\right)\left(\frac{\pi}{2}\cos\left(\frac{\pi}{2} t\right)+\sin\left(\frac{\pi}{2} t\right)\right)+(2 a-3 c)\left(\pi\cos\left(\pi t\right)+\sin\left(\pi t\right)\right)-
\\
\nonumber &-\left(\frac{1}{4}+2 a+b\right)\left(\frac{3\pi}{2}\cos\left(\frac{3\pi}{2} t\right)+\sin\left(\frac{3\pi}{2} t\right)\right)+
\\
\nonumber &+a\left(2\pi\cos\left(2\pi t\right)+\sin\left(2\pi t\right)\right)+b\left(\frac{5\pi}{2}\cos\left(\frac{5\pi}{2} t\right)+\sin\left(\frac{5\pi}{2} t\right)\right)+c\left(3\pi\cos\left(3\pi t\right)+  \sin\left(3\pi t\right)\right)
\end{align}
The corresponding values for the free parameters and for the cost functional are:
\begin{equation}
a\approx 1.0407,~ b\approx -0.312242,~ c\approx -0.0105136,~ \mathcal{C}_6 \approx 1.48099.
\end{equation}
\subsection{Exponential function family}
We consider in this paragraph that the trajectory can be written as follows:
\begin{equation}
\label{expeq}
x(t)=a \e^t+b \e^{-t}+c \e^{k t}+d \e^{-k t}
\end{equation}
where $a$, $b$, $c$, $d$ and $k$ are free parameters. From the boundary conditions $x(0)=0$, $x(1)=1$ and $\dot{x}(0)=\dot{x}(1)=0$, we have:
\begin{equation*}
B \vec{h}=\vec{s},
\end{equation*}
where
\begin{equation*}
B=
\left(
\begin{array}{cccc}
1 & 1 & 1 & 1
\\
\e & \e^{-1} & \e^k & \e^{-k}
\\
1 & -1 & k & -k
\\
\e & -\e^{-1} & k\e^k & -k\e^{-k}
\end{array}
\right),\quad
\vec{h}=
\left(
\begin{array}{cccc}
a
\\
b
\\
c
\\
d
\end{array}
\right),\quad
\vec{s}=
\left(
\begin{array}{cccc}
0
\\
1
\\
0
\\
0
\end{array}
\right)
\end{equation*}
The vector $\vec{h}$ is given by $\vec{h}=B^{-1} \vec{s}$, where
\begin{eqnarray*}
B^{-1}=\frac{1}{\det{B}}\left(
\begin{pmatrix}
\abs{B}_{11} & \abs{B}_{21} & \abs{B}_{31}  & \abs{B}_{41}
\\
\abs{B}_{12}  & \abs{B}_{22}  & \abs{B}_{32}  & \abs{B}_{42}
\\
\abs{B}_{13}  & \abs{B}_{23}  & \abs{B}_{33}  & \abs{B}_{43}
\\
\abs{B}_{14}  & \abs{B}_{24}  & \abs{B}_{34}  & \abs{B}_{44}
\end{pmatrix}
\right)
\end{eqnarray*}
where $\abs{B}_{11},~\abs{B}_{12},\dots$ are the cofactors of the $B$ matrix and $\det{B}$ is the determinant. We deduce that $\vec{h}$ can be expressed as:
\begin{equation*}
\vec{h}=\frac{1}{\det{B}}\left(
\begin{array}{cccc}
\abs{B}_{21}
\\
\abs{B}_{22}
\\
\abs{B}_{23}
\\
\abs{B}_{24}
\end{array}
\right)
\end{equation*}
Straightforward computations lead to the following determinant and cofactors:
\begin{equation*}
\begin{dcases}
\det{B}=-(1-k)^2 \e^{-(1+k)}-(1-k)^2 \e^{1+k}+(1+k)^2 \e^{-(1-k)}+(1+k)^2 \e^{1-k}-8 k
\\
\abs{B}_{21}=k\left(-2 \e^{-1}+(1+k) \e^{-k}+(1-k) \e^{k}\right)
\\
\abs{B}_{22}=k\left(-2 \e+(1-k) \e^{-k}+(1+k) \e^{k}\right)
\\
\abs{B}_{23}=k\left(\left(1+\frac{1}{k}\right) \e^{-1}+\left(1-\frac{1}{k}\right) \e-2\e^{-k}\right)
\\
\abs{B}_{24}=k\left(\left(1-\frac{1}{k}\right) \e^{-1}+\left(1+\frac{1}{k}\right) \e-2\e^{k}\right)
\end{dcases}
\end{equation*}
Whence the coefficients $a$, $b$, $c$ and $d$ read:
\begin{equation*}
\begin{dcases}
a=\frac{k}{D}\left(-2 \e^{-1}+(1+k) \e^{-k}+(1-k) \e^{k}\right)
\\
b=\frac{k}{D}\left(-2 \e+(1-k) \e^{-k}+(1+k) \e^{k}\right)
\\
c=\frac{k}{D}\left(\left(1+\frac{1}{k}\right) \e^{-1}+\left(1-\frac{1}{k}\right) \e-2\e^{-k}\right)
\\
d=\frac{k}{D}\left(\left(1-\frac{1}{k}\right) \e^{-1}+\left(1+\frac{1}{k}\right) \e-2\e^{k}\right)
\end{dcases}
\end{equation*}
with
$$
D=\det{B}=-(1-k)^2 \e^{-(1+k)}-(1-k)^2 \e^{1+k}+(1+k)^2 \e^{-(1-k)}+(1+k)^2 \e^{1-k}-8 k.
$$
The cost functional for $k=100$ is $\mathcal{C}=1.325271$.

\section{The regular optimal solution}
Pontryagin's equations in the regular case is given by:
\begin{equation*}
\dot{\vec{X}}=M_{\lambda} \vec{X},
\end{equation*}
where $\vec{X}=(y,z,p_y,p_z)^T$ is the state of the system, and $M_{\lambda}$ has the following form:
\begin{equation*}
M_{\lambda}=
\left(
\begin{array}{cccc}
-1 & 0 & 1/\lambda & 1/\lambda
\\
0 & 0 & 1/\lambda & 1/\lambda
\\
2 & -1 & 1 & 0
\\
-1 & 1 & 0 & 0
\end{array}
\right)
\end{equation*}
The analytic solution is:
\begin{equation*}
\vec{X}=\e^{M_{\lambda} t} \vec{X}(0),
\end{equation*}
where $\vec{X}(0)$ is the initial state of the system. The computation of the optimal solution requires the determination of the matrix exponential and of the initial adjoint state. We first compute  the eigenvalues and eigenvectors of the matrix $M_{\lambda}$:
\begin{equation*}
M_{\lambda}\vec{X}=\mu\vec{X}.
\end{equation*}
We obtain:
\begin{equation*}
\mu_1=-1,~ \mu_2=1,~ \mu_3=-\frac{1}{\sqrt{\lambda}},~ \mu_4=\frac{1}{\sqrt{\lambda}}
\end{equation*}
and
\begin{equation*}
X_1=c_1\left(
\begin{array}{c}
1
\\
0
\\
-1
\\
1
\end{array}
\right),
X_2=c_2\left(
\begin{array}{c}
1
\\
2
\\
2 \lambda-1
\\
1
\end{array}
\right)
\\
X_3=c_3\left(
\begin{array}{c}
1
\\
1-\sqrt{\lambda}
\\
-\sqrt{\lambda}
\\
\lambda
\end{array}
\right),
X_4=c_4\left(
\begin{array}{c}
1
\\
1+\sqrt{\lambda}
\\
\sqrt{\lambda}
\\
\lambda
\end{array}
\right)
\end{equation*}
where the $c_i$'s are normalization constants which can be set to 1. The eigenvector matrix $P$ can be written as:
\begin{equation*}
P=
\left(
\begin{array}{cccc}
1 & 1 & 1 & 1
\\
0 & 2 & 1-\sqrt{\lambda} & 1+\sqrt{\lambda}
\\
-1 & 2 \lambda-1 & -\sqrt{\lambda} & \sqrt{\lambda}
\\
1 & 1 & \lambda & \lambda
\end{array}
\right).
\end{equation*}
The matrix can be inverted if $\lambda\neq1$: $\det(P)=4 \lambda^{1/2} (1-\lambda)^2$. This is not a problem since we are interested in values of $\lambda$ close to $0$. We have:
\begin{equation*}
P^{-1}=
\begin{pmatrix}
\displaystyle\frac{1-2 \lambda}{2 (1-\lambda)} & -\displaystyle\frac{1}{2 (1-\lambda)} & \displaystyle\frac{1}{2 (1-\lambda)} & \displaystyle\frac{1}{1-\lambda}
\\
-\displaystyle\frac{1}{2 (1-\lambda)} & \displaystyle\frac{1}{2 (1-\lambda)} & -\displaystyle\frac{1}{2 (1-\lambda)} & 0
\\
\displaystyle\frac{1}{2 (1-\lambda)} & \displaystyle\frac{\sqrt{\lambda}}{2 (1-\lambda)} & -\displaystyle\frac{1}{2 \sqrt{\lambda} (1-\lambda)} & \displaystyle\frac{1}{2 (\lambda-\sqrt{\lambda})}
\\
\displaystyle\frac{1}{2 (1-\lambda)} & -\displaystyle\frac{\sqrt{\lambda}}{2 (1-\lambda)} & \displaystyle\frac{1}{2 \sqrt{\lambda} (1-\lambda)} & \displaystyle\frac{1}{2 (\lambda+\sqrt{\lambda})}
\end{pmatrix}
\end{equation*}
We obtain:
\begin{math}
M_{\lambda} P=P D,
\end{math}
where $D$ is the following diagonal matrix:
\begin{equation}
D=
\begin{pmatrix}
-1&0&0&0
\\
0&1&0&0
\\
0&0&-\displaystyle\frac{1}{\sqrt{\lambda}}&0
\\
0&0&0&\displaystyle\frac{1}{\sqrt{\lambda}}.
\end{pmatrix},
\end{equation}
It is then straightforward to get $N(t)$:
\begin{equation}
N(t)=e^{M_{\lambda} t}=\left(N_1(t)~ N_2(t)~ N_3(t)~ N_4(t)\right),
\end{equation}
where the vectors $N_1(t)$, $N_2(t)$, $N_3(t)$ and $N_4(t)$ are defined as:
\begin{eqnarray}
\nonumber N_1(t)&=&\frac{1}{2 (1-\lambda)}
\begin{pmatrix}
\displaystyle (1-2 \lambda) \e^{-t}-\e^{t}+\e^{-\frac{t}{\sqrt{\lambda}}}+\e^{\frac{t}{\sqrt{\lambda}}}
\\[0.5cm]
\displaystyle -\e^{-t}+\e^{t}+\sqrt{\lambda} \e^{-\frac{t}{\sqrt{\lambda}}}-\sqrt{\lambda} \e^{\frac{t}{\sqrt{\lambda}}}
\\[0.5cm]
\displaystyle\e^{-t}-\e^{t}-\frac{1}{\sqrt{\lambda}} \e^{-\frac{t}{\sqrt{\lambda}}}+\frac{1}{\sqrt{\lambda}} \e^{\frac{t}{\sqrt{\lambda}}}
\\[0.5cm]
\displaystyle 2 \e^{-t}-\left(1+\frac{1}{\sqrt{\lambda}}\right) \e^{-\frac{t}{\sqrt{\lambda}}}-\left(1-\frac{1}{\sqrt{\lambda}}\right) \e^{\frac{t}{\sqrt{\lambda}}}
\end{pmatrix}
\\
\nonumber N_2(t)&=&\frac{1}{2 (1-\lambda)}
\begin{pmatrix}
\displaystyle -2 \e^{t}+\left(1-\sqrt{\lambda}\right) \e^{-\frac{t}{\sqrt{\lambda}}}+\left(1+\sqrt{\lambda}\right) \e^{\frac{t}{\sqrt{\lambda}}}
\\[0.5cm]
\displaystyle 2 \e^{t}+\left(\sqrt{\lambda}-\lambda\right) \e^{-\frac{t}{\sqrt{\lambda}}}-\left(\sqrt{\lambda}+\lambda\right) \e^{\frac{t}{\sqrt{\lambda}}}
\\[0.5cm]
\displaystyle -2 \e^{t}+\left(1-\frac{1}{\sqrt{\lambda}}\right) \e^{-\frac{t}{\sqrt{\lambda}}}+\left(1+\frac{1}{\sqrt{\lambda}}\right) \e^{\frac{t}{\sqrt{\lambda}}}
\\[0.5cm]
\displaystyle\left(\frac{1}{\sqrt{\lambda}}-\sqrt{\lambda}\right) \left(-\e^{-\frac{t}{\sqrt{\lambda}}}+\e^{\frac{t}{\sqrt{\lambda}}}\right)
\end{pmatrix}
\\
\nonumber N_3(t)&=&\frac{1}{2 (1-\lambda)}
\begin{pmatrix}
\displaystyle -(1-2 \lambda)\e^{-t}+(1-2 \lambda)\e^{t}-\sqrt{\lambda}\e^{-\frac{t}{\sqrt{\lambda}}}+\sqrt{\lambda}\e^{\frac{t}{\sqrt{\lambda}}}
\\[0.5cm]
\displaystyle \e^{-t}-(1-2 \lambda)\e^{t}-\lambda\e^{-\frac{t}{\sqrt{\lambda}}}-\lambda\e^{\frac{t}{\sqrt{\lambda}}}
\\[0.5cm]
\displaystyle -\e^{-t}+(1-2 \lambda)\e^{t}+\e^{-\frac{t}{\sqrt{\lambda}}}+\e^{\frac{t}{\sqrt{\lambda}}}
\\[0.5cm]
\displaystyle -2 \e^{-t}+\left(1+\sqrt{\lambda}\right)\e^{-\frac{t}{\sqrt{\lambda}}}+\left(1-\sqrt{\lambda}\right)\e^{\frac{t}{\sqrt{\lambda}}}
\end{pmatrix}
\\
\nonumber N_4(t)&=&\frac{1}{2 (1-\lambda)}
\begin{pmatrix}
\displaystyle (1-2 \lambda)\e^{-t}-\e^{t}+\lambda\e^{-\frac{t}{\sqrt{\lambda}}}+\lambda\e^{\frac{t}{\sqrt{\lambda}}}
\\[0.5cm]
\displaystyle -\e^{-t}+\e^{t}+\lambda\sqrt{\lambda}\e^{-\frac{t}{\sqrt{\lambda}}}-\lambda\sqrt{\lambda}\e^{\frac{t}{\sqrt{\lambda}}}
\\[0.5cm]
\displaystyle \e^{-t}-\e^{t}-\sqrt{\lambda}\e^{-\frac{t}{\sqrt{\lambda}}}+\sqrt{\lambda}\e^{\frac{t}{\sqrt{\lambda}}}
\\[0.5cm]
\displaystyle 2 \e^{-t}-\left(\sqrt{\lambda}+\lambda\right)\e^{-\frac{t}{\sqrt{\lambda}}}+\left(\sqrt{\lambda}-\lambda\right)\e^{\frac{t}{\sqrt{\lambda}}}
\end{pmatrix}
\end{eqnarray}
From Pontryagin's equations, it is also easy to deduce that:
\begin{equation*}
\begin{pmatrix}
\displaystyle N_{13}(t) & N_{14}(t)
\\
\displaystyle  N_{23}(t) & N_{24}(t)
\end{pmatrix}
\begin{pmatrix}
\displaystyle  p_y(0)
\\
\displaystyle  p_z(0)
\end{pmatrix}
=
\begin{pmatrix}
\displaystyle  y(t)-N_{11}(t) y(0)-N_{12}(t) z(0)
\\
\displaystyle  z(t)-N_{21}(t) y(0)-N_{22}(t)z(0)
\end{pmatrix}
\end{equation*}
where the $N_{ij}$s are the matrix elements of $N$. Using the boundary conditions $y(0)=z(0)=0$, we arrive at:
\begin{equation*}
\label{mat_eq_2}
\begin{pmatrix}
\displaystyle  y(t)
\\
\displaystyle  z(t)
\end{pmatrix}
=
\begin{pmatrix}
\displaystyle N_{13}(t) & N_{14}(t)
\\
\displaystyle  N_{23}(t) & N_{24}(t)
\end{pmatrix}
\begin{pmatrix}
\displaystyle  p_y(0)
\\
\displaystyle  p_z(0)
\end{pmatrix}.
\end{equation*}
Similarly, we get:
\begin{equation*}
\label{mat_eq_3}
\begin{pmatrix}
\displaystyle  p_y(t)
\\
\displaystyle  p_z(t)
\end{pmatrix}
=
\begin{pmatrix}
\displaystyle N_{33}(t) & N_{34}(t)
\\
\displaystyle  N_{43}(t) & N_{44}(t)
\end{pmatrix}
\begin{pmatrix}
\displaystyle  p_y(0)
\\
\displaystyle  p_z(0)
\end{pmatrix}.
\end{equation*}
We have established all the intermediate results to find the initial adjoint state. From $y(T)=0$, and $z(T)=1$, we have:
\begin{equation*}
\begin{pmatrix}
\displaystyle  p_y(0)
\\
\displaystyle  p_z(0)
\end{pmatrix}
=
A^{-1}
\begin{pmatrix}
\displaystyle  0
\\
\displaystyle  1
\end{pmatrix},
\end{equation*}
where
\begin{equation*}
A=
\begin{pmatrix}
\displaystyle N_{13}(T) & N_{14}(T)
\\
\displaystyle  N_{23}(T) & N_{24}(T)
\end{pmatrix},
\end{equation*}
It can be shown that:
\begin{equation*}
A^{-1}=
\alpha
\begin{pmatrix}
\displaystyle\left(\frac{1}{\sqrt{\lambda}}-\sqrt{\lambda}\right) \left(-\e^{-\frac{T}{\sqrt{\lambda}}}+\e^{\frac{T}{\sqrt{\lambda}}}\right) &
-\displaystyle 2 \e^{-T}+(1+\frac{1}{\sqrt{\lambda}}) \e^{-\frac{T}{\sqrt{\lambda}}}+(1-\frac{1}{\sqrt{\lambda}}) \e^{\frac{T}{\sqrt{\lambda}}}
\\[0.5cm]
\displaystyle 2 \e^{T}-\left(1-\frac{1}{\sqrt{\lambda}}\right) \e^{-\frac{T}{\sqrt{\lambda}}}-\left(1+\frac{1}{\sqrt{\lambda}}\right) \e^{\frac{T}{\sqrt{\lambda}}} &
\displaystyle\e^{-T}-\e^{T}-\frac{1}{\sqrt{\lambda}} \e^{-\frac{T}{\sqrt{\lambda}}}+\frac{1}{\sqrt{\lambda}} \e^{\frac{T}{\sqrt{\lambda}}}
\end{pmatrix},
\end{equation*}
with $\alpha=\displaystyle \frac{1}{2 (1-\lambda) \det(A)}$ and
\begin{eqnarray*}
& & \det(A)=\frac{1}{4 (1-\lambda)^2 \sqrt{\lambda}}\left[\left(1-\sqrt{\lambda}\right)^2 \e^{-\left(1+\frac{1}{\sqrt{\lambda}}\right) T}-\left(1+\sqrt{\lambda}\right)^2 \e^{-\left(1-\frac{1}{\sqrt{\lambda}}\right) T}\right.
\\
& & -\left.\left(1+\sqrt{\lambda}\right)^2 \e^{\left(1-\frac{1}{\sqrt{\lambda}}\right) T}+\left(1-\sqrt{\lambda}\right)^2 \e^{\left(1+\frac{1}{\sqrt{\lambda}}\right) T}+8\sqrt{\lambda}\right].
\end{eqnarray*}
We finally obtain the initial adjoint state:
\begin{equation*}
\begin{pmatrix}
p_y(0)
\\
p_z(0)
\end{pmatrix}
=\frac{1}{2 (1-\lambda) det(A)}
\begin{pmatrix}
-\displaystyle 2 \e^{-T}+\left(1+\frac{1}{\sqrt{\lambda}}\right) \e^{-\frac{T}{\sqrt{\lambda}}}+\left(1-\frac{1}{\sqrt{\lambda}}\right) \e^{\frac{T}{\sqrt{\lambda}}}
\\[0.5cm]
\displaystyle\e^{-T}-\e^{T}-\frac{1}{\sqrt{\lambda}} \e^{-\frac{T}{\sqrt{\lambda}}}+\frac{1}{\sqrt{\lambda}} \e^{\frac{T}{\sqrt{\lambda}}}
\end{pmatrix}.
\end{equation*}
It is then straightforward to derive the dynamics of the system:
\begin{equation}\label{state}
\begin{pmatrix}
y
\\
z
\\
x
\\
p_y
\\
p_z
\\
v_R
\end{pmatrix}
=\a
\begin{pmatrix}
y_1 & y_2 & y_3 & y_4
\\
z_1 & z_2 & z_3 & z_4
\\
x_1 & x_2 & x_3 & x_4
\\
p_{y_1} & p_{y_2} & p_{y_3} & p_{y_4}
\\
p_{z_1} & p_{z_2} & p_{z_3} & p_{z_4}
\\
v_1 & v_2 & v_3 & v_4
\end{pmatrix}
\begin{pmatrix}
\e^{-t}
\\
\e^{t}
\\
\e^{-\frac{t}{\sqrt{\lambda}}}
\\
\e^{\frac{t}{\sqrt{\lambda}}}
\end{pmatrix}
\end{equation}
where
\begin{equation}
\label{gotha}
\a=\frac{1}{4 (1-\lambda)^2 \textrm{det}(A)}.
\end{equation}
It is also worth to recall that:
\begin{gather*}
x=z-y
\\
v_R=\frac{p_y+p_z}{\lambda}
\end{gather*}
The $y_1$, $y_2$, $y_3$, $y_4$ coefficients can be expressed as follows:
\begin{eqnarray}
\nonumber y_1&=&-2\e^{T}+\left(1-\frac{1}{\sqrt{\lambda}}\right)\e^{-\frac{T}{\sqrt{\lambda}}}+\left(1+\frac{1}{\sqrt{\lambda}}\right)\e^{\frac{T}{\sqrt{\lambda}}}
\\
\nonumber y_2&=&2\e^{-T}-\left(1+\frac{1}{\sqrt{\lambda}}\right)\e^{-\frac{T}{\sqrt{\lambda}}}-\left(1-\frac{1}{\sqrt{\lambda}}\right)\e^{\frac{T}{\sqrt{\lambda}}}
\\
\nonumber y_3&=&-\left(1-\frac{1}{\sqrt{\lambda}}\right)\e^{-T}+\left(1+\frac{1}{\sqrt{\lambda}}\right)\e^{T}-\frac{2}{\sqrt{\lambda}}\e^{\frac{T}{\sqrt{\lambda}}}
\\
\nonumber y_4&=&-\left(1+\frac{1}{\sqrt{\lambda}}\right)\e^{-T}+\left(1-\frac{1}{\sqrt{\lambda}}\right)\e^{T}+\frac{2}{\sqrt{\lambda}}\e^{-\frac{T}{\sqrt{\lambda}}}
\end{eqnarray}
All the other coefficients $z_1$, $z_2$, $\dots$, $v_4$ can be expressed in terms of $y_1$, $y_2$, $y_3$ and $y_4$:
\begin{equation*}
\label{multipliers}
\left[
\begin{array}{ccccc}
\displaystyle z_1 & \displaystyle z_2 & \displaystyle z_3 & \displaystyle z_4
\\
\displaystyle x_1 & \displaystyle x_2 & \displaystyle x_3 & \displaystyle x_4
\\
\displaystyle p_{y1} & \displaystyle p_{y2} & \displaystyle p_{y3} & \displaystyle p_{y4}
\\
\displaystyle p_{z1} & \displaystyle p_{z2} & \displaystyle p_{z3} & \displaystyle p_{z4}
\\
\displaystyle v_1 & \displaystyle v_2 & \displaystyle v_3 & \displaystyle v_4
\end{array}
\right]
=
\left[
\begin{array}{ccccc}
\displaystyle 0 & \displaystyle 2 y_2 & \displaystyle \left(1-\sqrt{\lambda}\right) y_3 & \displaystyle \left(1+\sqrt{\lambda}\right) y_4
\\
\displaystyle -y_1 & \displaystyle y_2 & \displaystyle -\sqrt{\lambda} y_3 & \displaystyle \sqrt{\lambda} y_4
\\
\displaystyle -y_1 & \displaystyle -\left(1-2 \lambda\right) y_2 & \displaystyle -\sqrt{\lambda} y_3 & \displaystyle \sqrt{\lambda} y_4
\\
\displaystyle y_1 & \displaystyle y_2 & \displaystyle \lambda y_3 & \displaystyle \lambda y_4
\\
\displaystyle 0 & \displaystyle 2 y_2 & \displaystyle \left(1-\frac{1}{\sqrt{\lambda}}\right) y_3 & \displaystyle \left(1+\frac{1}{\sqrt{\lambda}}\right) y_4
\end{array}
\right]
\end{equation*}

\noindent \textbf{The regularized cost functional}\\
Having established the optimal regular solution, the next step is to compute the corresponding cost functional. We have:
\begin{equation}
\mathcal{C}_R=\int\limits_0^T\!\!G dt =\int\limits_0^T\!\!\left[(z-y)^2+y^2+\lambda v_R^2\right]\! dt
\end{equation}
The different terms can be expressed analytically:
\begin{equation*}
\lambda v_R^2=\a^2\left(w_2^2\e^{2 t}+w_3^2\e^{-\frac{2 t}{\sqrt{\lambda}}}+w_4^2\e^{\frac{2 t}{\sqrt{\lambda}}}+2\left[w_2 w_3\e^{\left(1-\frac{1}{\sqrt{\lambda}}\right) t}+w_2 w_4\e^{\left(1+\frac{1}{\sqrt{\lambda}}\right) t}+w_3 w_4\right]\right),
\end{equation*}
where
\begin{equation*}
w_2=2\sqrt{\lambda} y_2,~ w_3=-\left(1-\sqrt{\lambda}\right) y_3, ~ w_4=\left(1+\sqrt{\lambda}\right) y_4
\end{equation*}
We also have:
\begin{eqnarray*}
& & y^2=\a^2(y_1^2\e^{-2 t}+y_2^2\e^{2 t}+y_3^2\e^{-\frac{2 t}{\sqrt{\lambda}}}+y_4^2\e^{\frac{2 t}{\sqrt{\lambda}}}+2(y_1 y_2+y_1 y_3\e^{-(1+\frac{1}{\sqrt{\lambda}})t}\\
& & +y_1y_4\e^{-(1-\frac{1}{\sqrt{\lambda}}) t} +y_2 y_3\e^{(1-\frac{1}{\sqrt{\lambda}}) t}+y_2 y_4\e^{(1+\frac{1}{\sqrt{\lambda}}) t}+y_3 y_4))
\end{eqnarray*}
\begin{eqnarray*}
& & (z-y)^2=\a^2(x_1^2\e^{-2 t}+x_2^2\e^{2 t}+x_3^2\e^{-\frac{2 t}{\sqrt{\lambda}}}+x_4^2\e^{\frac{2 t}{\sqrt{\lambda}}}+2(x_1 x_2+x_1 x_3\e^{-(1+\frac{1}{\sqrt{\lambda}})t}\\
& & +x_1 x_4\e^{-(1-\frac{1}{\sqrt{\lambda}}) t}+x_2 x_3\e^{(1-\frac{1}{\sqrt{\lambda}})t}+x_2 x_4\e^{(1+\frac{1}{\sqrt{\lambda}}) t}+x_3 x_4 ))
\end{eqnarray*}
Simple manipulations yield to:
\begin{eqnarray*}
& & G=\a^2( G_1 \e^{-2 t}+G_2 \e^{2 t}+G_3 \e^{-\frac{2 t}{\sqrt{\lambda}}}+G_4 \e^{\frac{2 t}{\sqrt{\lambda}}}+2(G_{13} \e^{-(1+\frac{1}{\sqrt{\lambda}})t}\\
& & +G_{14} \e^{-(1-\frac{1}{\sqrt{\lambda}})t}+G_{23} \e^{(1-\frac{1}{\sqrt{\lambda}})t}
+G_{24} \e^{(1+\frac{1}{\sqrt{\lambda}})t})),
\end{eqnarray*}
where
\begin{align}
\nonumber & G_1=2 y_1^2,\quad G_2=2\left(1+2\lambda\right) y_2^2,\quad G_3=2\left(1-\sqrt{\lambda}+\lambda\right) y_3^2,\quad G_4=2\left(1+\sqrt{\lambda}+\lambda\right) y_4^2
\\
\nonumber & G_{13}=\left(1+\sqrt{\lambda}\right) y_1 y_3,\quad
G_{14}=\left(1-\sqrt{\lambda}\right) y_1 y_4
\\
\nonumber & G_{23}=\left(1-\sqrt{\lambda}\right)\left(1-2\sqrt{\lambda}\right) y_2 y_3,
\quad G_{24}=\left(1+\sqrt{\lambda}\right)\left(1+2\sqrt{\lambda}\right) y_2 y_4.
\end{align}
A simple integration leads to:
\begin{align}
& \nonumber \mathcal{C}_R=\a^2\left\{\frac{1}{2}\left[-\left(\e^{-2 T}-1\right) G_1+\left(\e^{2 T}-1\right) G_2\right]-\frac{\sqrt{\lambda}}{2}\left[\left(\e^{-\frac{2 T}{\sqrt{\lambda}}}-1\right) G_3-\left(\e^{\frac{2 T}{\sqrt{\lambda}}}-1\right) G_4\right]-\right.
\\
& \nonumber \left.-2\frac{\sqrt{\lambda}}{1+\sqrt{\lambda}}\left[\left(\e^{-\left(1+\frac{1}{\sqrt{\lambda}}\right) T}-1\right) G_{13}-\left(\e^{\left(1+\frac{1}{\sqrt{\lambda}}\right) T}-1\right) G_{24}\right]+\right.
\\
& \nonumber \left.+2\frac{\sqrt{\lambda}}{1-\sqrt{\lambda}}\left[\left(\e^{-\left(1-\frac{1}{\sqrt{\lambda}}\right) T}-1\right) G_{14}-\left(\e^{\left(1-\frac{1}{\sqrt{\lambda}}\right) T}-1\right) G_{23}\right]\right\}
\end{align}

\section{Comparison between the STA exponential and the optimal regular solutions}
We study in this paragraph the correspondence between the STA exponential and the optimal regular solutions. For the regular optimal trajectory, we have:
\begin{equation}
\label{regcom}
x_R(t)=\a \left(x_1 \e^{-t}+x_2 \e^{t}+x_3 \e^{-\frac{t}{\sqrt{\lambda}}}+x_4 \e^{\frac{t}{\sqrt{\lambda}}}\right)
\end{equation}
where the coefficients $x_1$, $x_2$, $x_3$, $x_4$ and $\a$ are given by:
\begin{gather}
\nonumber x_1=-y_1=2 \e-\left(1-\frac{1}{\sqrt{\lambda}}\right)\e^{-\frac{1}{\sqrt{\lambda}}}-\left(1+\frac{1}{\sqrt{\lambda}}\right)\e^{\frac{1}{\sqrt{\lambda}}}
\\
\nonumber x_2=y_2=2 \e^{-1}-\left(1+\frac{1}{\sqrt{\lambda}}\right)\e^{-\frac{1}{\sqrt{\lambda}}}-\left(1-\frac{1}{\sqrt{\lambda}}\right)\e^{\frac{1}{\sqrt{\lambda}}}
\\
\nonumber x_3=-\sqrt{\lambda} y_3=-\left(1-\sqrt{\lambda}\right)\e^{-1}-\left(1+\sqrt{\lambda}\right)\e+2\e^{\frac{1}{\sqrt{\lambda}}}
\\
\nonumber x_4=\sqrt{\lambda} y_4=-\left(1+\sqrt{\lambda}\right)\e^{-1}-\left(1-\sqrt{\lambda}\right)\e+2\e^{-\frac{1}{\sqrt{\lambda}}}
\\
\nonumber \a=\frac{\sqrt{\lambda}}{\left(1-\sqrt{\lambda}\right)^2 \e^{-\left(1+\frac{1}{\sqrt{\lambda}}\right)}-\left(1+\sqrt{\lambda}\right)^2 \e^{-\left(1-\frac{1}{\sqrt{\lambda}}\right)}
-\left(1+\sqrt{\lambda}\right)^2 \e^{1-\frac{1}{\sqrt{\lambda}}}+\left(1-\sqrt{\lambda}\right)^2 \e^{1+\frac{1}{\sqrt{\lambda}}}+8\sqrt{\lambda}}
\end{gather}
When $k\equiv \frac{1}{\sqrt{\lambda}}$, a relation can be established between $\textrm{det}{B}$ and $\a$:
\begin{align}
\nonumber \det{B}=-\frac{1}{\lambda}\left[(1-\sqrt{\lambda})^2 \e^{-\left(1+\frac{1}{\sqrt{\lambda}}\right)}-(1+\sqrt{\lambda})^2 \e^{-\left(1-\frac{1}{\sqrt{\lambda}}\right)}-\right.
\\
\nonumber \left.-(1+\sqrt{\lambda})^2 \e^{\left(1-\frac{1}{\sqrt{\lambda}}\right)}+(1-\sqrt{\lambda})^2 \e^{\left(1+\frac{1}{\sqrt{\lambda}}\right)}+8\sqrt{\lambda}\right],
\end{align}
or
\begin{equation*}
\frac{1}{\det{B}}=-\sqrt{\lambda}\a
\end{equation*}
the $a$, $b$, $c$ and $d$ coefficients can be expressed as:
\begin{equation*}
\begin{dcases}
b=-\frac{1}{\det{B}}\frac{1}{\sqrt{\lambda}} x_1=\a x_1
\\
a=-\frac{1}{\det{B}}\frac{1}{\sqrt{\lambda}} x_2=\a x_2
\\
d=-\frac{1}{\det{B}}\frac{1}{\sqrt{\lambda}} x_3=\a x_3
\\
c=-\frac{1}{\det{B}}\frac{1}{\sqrt{\lambda}} x_4=\a x_4
\end{dcases}
\end{equation*}
Consequently the STA exponential solution reduces to:
\begin{equation}
\label{expcom}
x(t)=\a \left(x_1 \e^{-t}+x_2 \e^{t}+x_3 \e^{-\frac{t}{\sqrt{\lambda}}}+x_4 \e^{\frac{t}{\sqrt{\lambda}}}\right)
\end{equation}
From Eq.~\eqref{regcom} and Eq.~\eqref{expcom}, we deduce that the regular and exponential solutions coincide.

\end{document}